\documentclass[english]{article}
\usepackage[T1]{fontenc}
\usepackage[latin9]{inputenc}
\usepackage{graphicx}

\makeatletter

\providecommand{\tabularnewline}{\\}

\makeatother

\makeatother

\usepackage{babel}

\makeatother

\usepackage{babel}

\begin{document}

\title{3D Volume Calculation For the Marching Cubes Algorithm in Cartesian
Coordinates}

\author{Shuqiang Wang}

\date{06/20/2013}

\maketitle
\tableofcontents{}
\begin{abstract}
From a scalar field defined at the corner of a cube, an isosurface
can be extracted using the Marching Cube algorithm. The isosurface
separates the cell into two or more partial cells. A similar situation
arises when an material interface in the Front Tracking method cuts
through the computational cells. A popular method to calculate the
volumes of the partial cells is to first partition the cells into
tetrahedra and then sum together the volumes of the tetrahedra for
the corresponding partial cells. In this paper, the divergence theorem
is used to calculate the volumes of the partial cells generated by
the Marching Cubes algorithm. This method is both more robust and
efficient compared with the tetrahedralization approach.
\end{abstract}

\section{Introduction}

The Marching Cubes algorithm \cite{LorCli87} was developed to reconstruct
the interface using the volumetric data. The isosurface separates
the cell into two or more partial cells. A similar situation arises
when an material interface in the Front Tracking method cuts through
the computational cells. Most of the publications are concerning on
how to deal with interface reconstructions, while only few papers
exists to show how to calculate the partial volumes enclosed by the
interface. Current volume calculation methods \cite{KimPestieauGlimm2007,HareGroshSchmitt1999}
often first tetrahedralize the partial volumes and then calculate
the volume using the tetrahedra. This approach is difficult to write
computer code to deal with all possible cases. An alternative method
to the Marching Cubes algorithm is the generalized Marching Cubes
method \cite{Bloomenthal1997} which first partitions the cubic cell
into tetrahedra, and then uses the generalized Marching Cubes method
on the tetrahedra \cite{Bloomenthal1997} to recover the interface,
leading to few cases compared with the Marching Cubes algorithm. However,
the disadvantage of using this approach to calculate the partial volumes
is that we need crossing information on the cell face diagonals and
the main cube diagonal which might be difficult to obtain. In this
paper, we show how to calcalute the cell partial volumes by using
the divergence theorm. This approach seems to be much more robust
and computing efficient than the other approach. It is also very easy
to be coded.

As applications of this algorithm, we can/have used it in our embedded
boundary method (a finite volume method) coupled with the Front Tracking
method for solving the elliptic/parabolic interface problem \cite{WangSamulyakGuo2010},
two-phase incompressible flow \cite{WangGlimmSamulyakJiaoDiao2013},
and magnetohydrodynamic flows \cite{GuoWangSamulyak2013}.

\section{Method}

The general idea of using the divergence theorem to calculate the
volume is very simple. The divergence theorem is \begin{equation}
\int_{\Omega}\nabla\cdot\overrightarrow{f}dv=\oint_{\partial\Omega}\overrightarrow{n}\cdot\overrightarrow{f}ds.\label{eq:Divergence-Theorem}\end{equation}
If we can find a vector $\overrightarrow{f}$ such that \[
\nabla\cdot\overrightarrow{f}=a\]
where $a$ is a constant, then we have \[
\int_{\Omega}adv=\oint_{\partial\Omega}\overrightarrow{n}\cdot\overrightarrow{f}ds\]
or \[
Volume(\Omega)=\int_{\Omega}dv=\frac{\oint_{\partial\Omega}\overrightarrow{n}\cdot\overrightarrow{f}ds}{a}\]

\subsection{Volume Calculations in Cartesian Coordinate}

In a Cartesian coordinate, the divergence operator is defined as \[
\nabla\cdot\overrightarrow{f}=\frac{\partial f_{x}}{\partial x}+\frac{\partial f_{y}}{\partial y}+\frac{\partial f_{z}}{\partial z}\]
 where $\overrightarrow{f}=(f_{x},f_{y},f_{z})^{T}$.

To calculate the volume of the domain $\Omega$, we can let $\overrightarrow{f}=(x,y,z)^{T}$
(which is not unique) and then we have \begin{eqnarray*}
\int_{\Omega}\frac{\partial x}{\partial x}+\frac{\partial y}{\partial y}+\frac{\partial z}{\partial z}dv & = & \oint_{\partial\Omega}\overrightarrow{n}\cdot\overrightarrow{f}dS\\
\int_{\Omega}3dv & = & \oint_{\partial\Omega}\overrightarrow{n}\cdot\overrightarrow{f}dS\end{eqnarray*}
 Therefore, we can calculate the volume using the surface integration\begin{equation}
Volume(\Omega)=\int_{\Omega}dv=\frac{\int_{\partial\Omega}\overrightarrow{n}\cdot\overrightarrow{f}dS}{3}.\label{eq:volume-calculation-for-Cartesian-Coordinate}\end{equation}

Thus, if the domain boundary $\partial\Omega$ consists of triangle
meshs in the cartesian coordinate, we only need to use second order
accuarate quadrature rule to obtain the exact volume.

\subsection{Templates and Coding}

In this paper, we use the $23$ unique cube configurations (or cases)
in \cite{NewmanYi2006} as the $15$ cube configurations in the original
Marching Cubes algorithm \cite{LorCli87} has consistency issue \cite{NewmanYi2006}.

\begin{figure}
\caption{\label{fig:Case00}Case00}

\centering{}\includegraphics[scale=0.3]{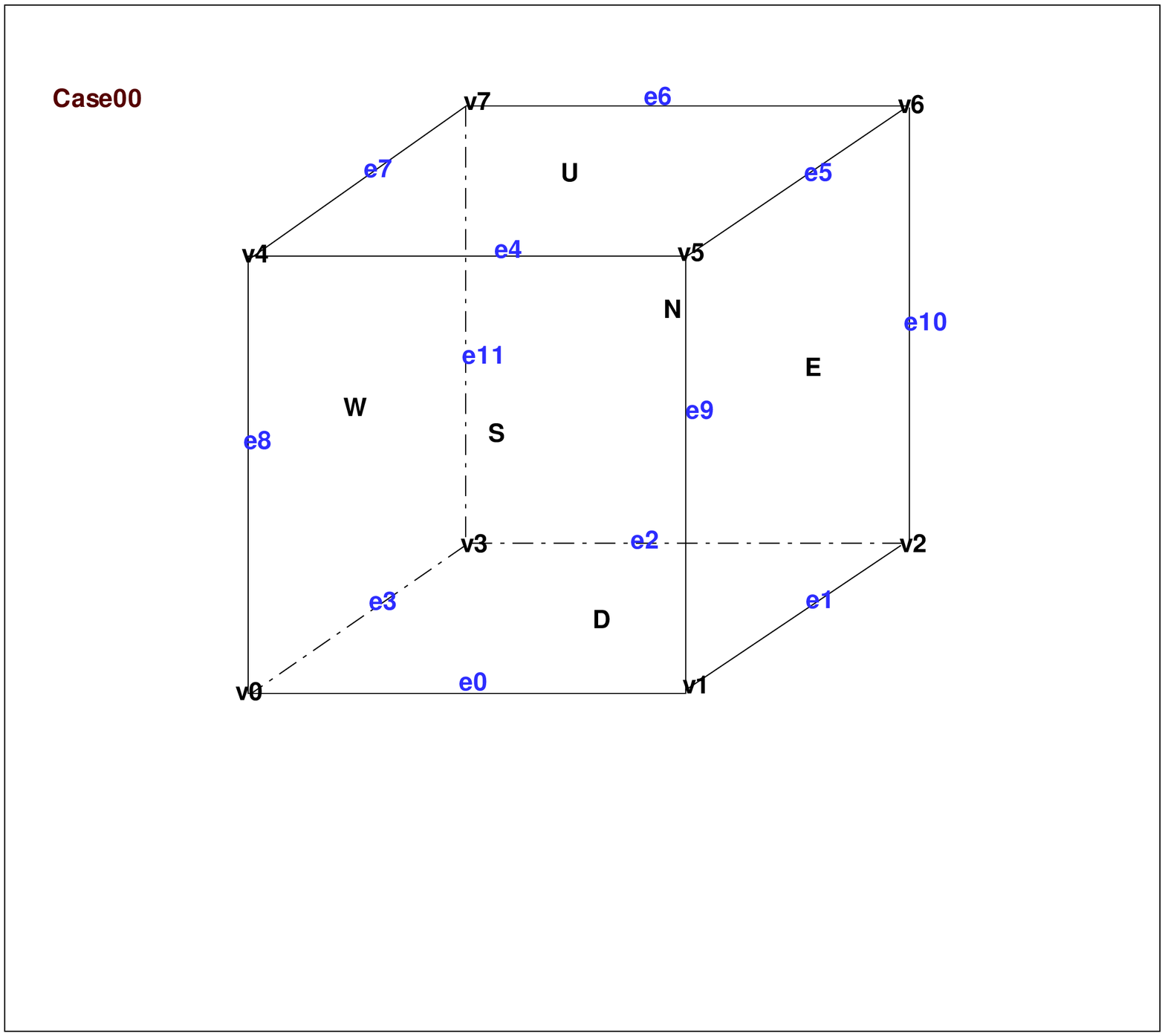} 
\end{figure}

Figure \ref{fig:Case00} shows the vertex, edge and face labeling
scheme which is different from \cite{NewmanYi2006}. Note that we
use a two letter word of pattern $v\bullet$ to denote the vertex
labeling and a two letter word of pattern $e\bullet$ to denote the
edge labeling. This kind of labeling makes it easy to write the code
to calculate the volume of the cubes since we can define those two
letter words as some constants and use them as the indices of the
coordinate arrays. The six faces of the cube are labeled as $W$ for
West (with vertex $v0,$$v3$, $v7$, $v4$), $E$ for East (with
vertex $v1,$$v2$, $v6$, $v5$), $S$ for South (with vertex $v0,$
$v1$, $v5$, $v4$), $N$ for North (with vertex $v3,$$v2$, $v6$,
$v7$), $D$ for Down (with vertex $v0,$$v1$, $v2$, $v3$), and
$U$ for Upper (with vertex $v4,$$v5$, $v6$, $v7$). On each vertex,
the componet is either $0$ or $1$. For component $1$ vertex, we
draw a circle on the vertex.

To calculate the partial volumes for a cube with two different components,
we calculate the volume for one of the components first, say $Volume_{1}$.
Then we use the total volume of the cube to calculate the volume for
the other component: \[
Volume_{0}=Volume_{total}-Volume_{1}\]

In order to calculate the volume using the surface integration in
equation (\ref{eq:volume-calculation-for-Cartesian-Coordinate}),
we need to find a closed surface enclosing that component and the
triangulation of the closed surface

\begin{figure}
\caption{\label{fig:Case01-Case02}Case01 and Case02}

\includegraphics[scale=0.3]{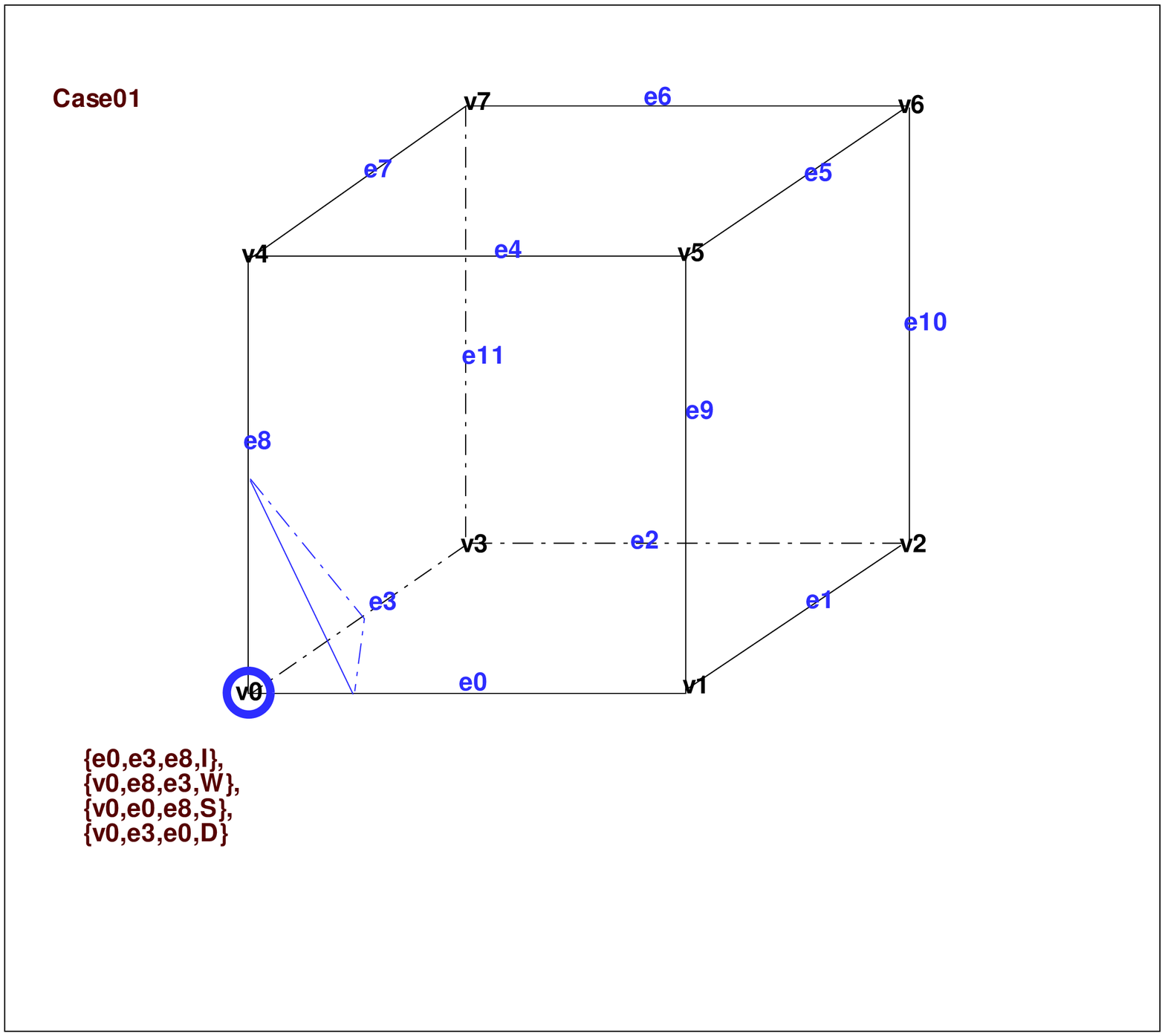} \includegraphics[scale=0.3]{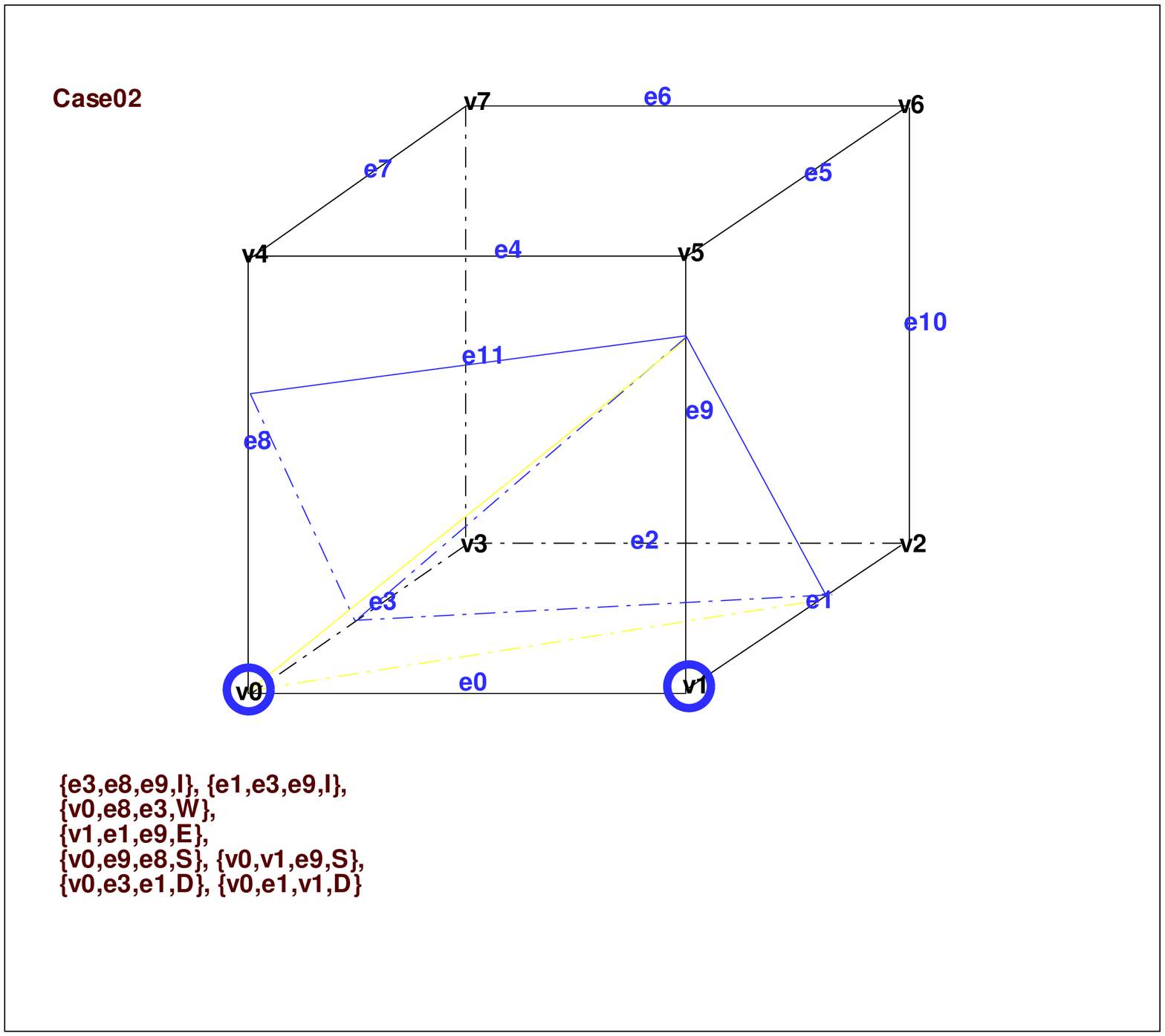} 
\end{figure}

\begin{figure}
\caption{\label{fig:Case03-Case04}Case03 and Case04}

\includegraphics[scale=0.3]{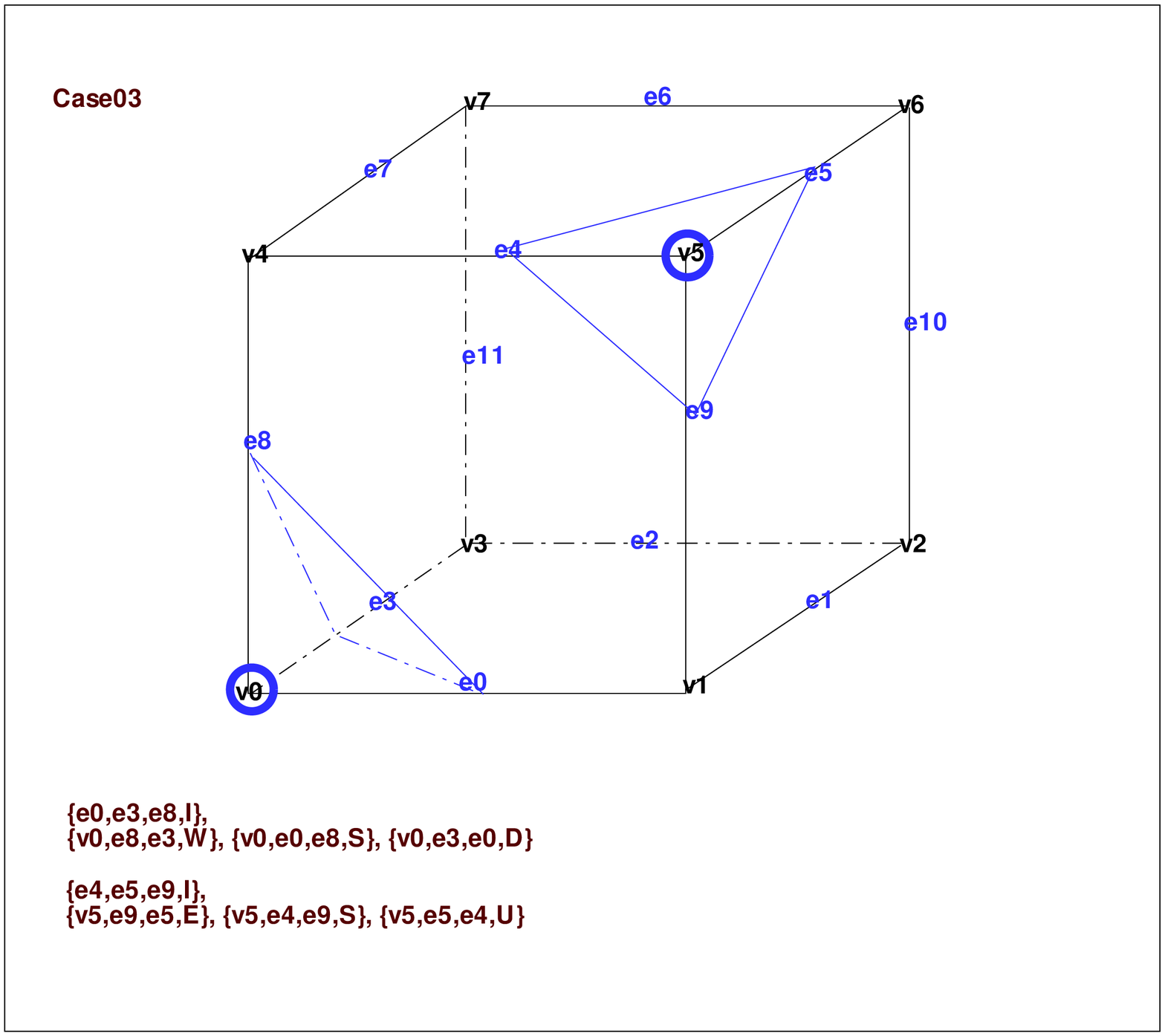} \includegraphics[scale=0.3]{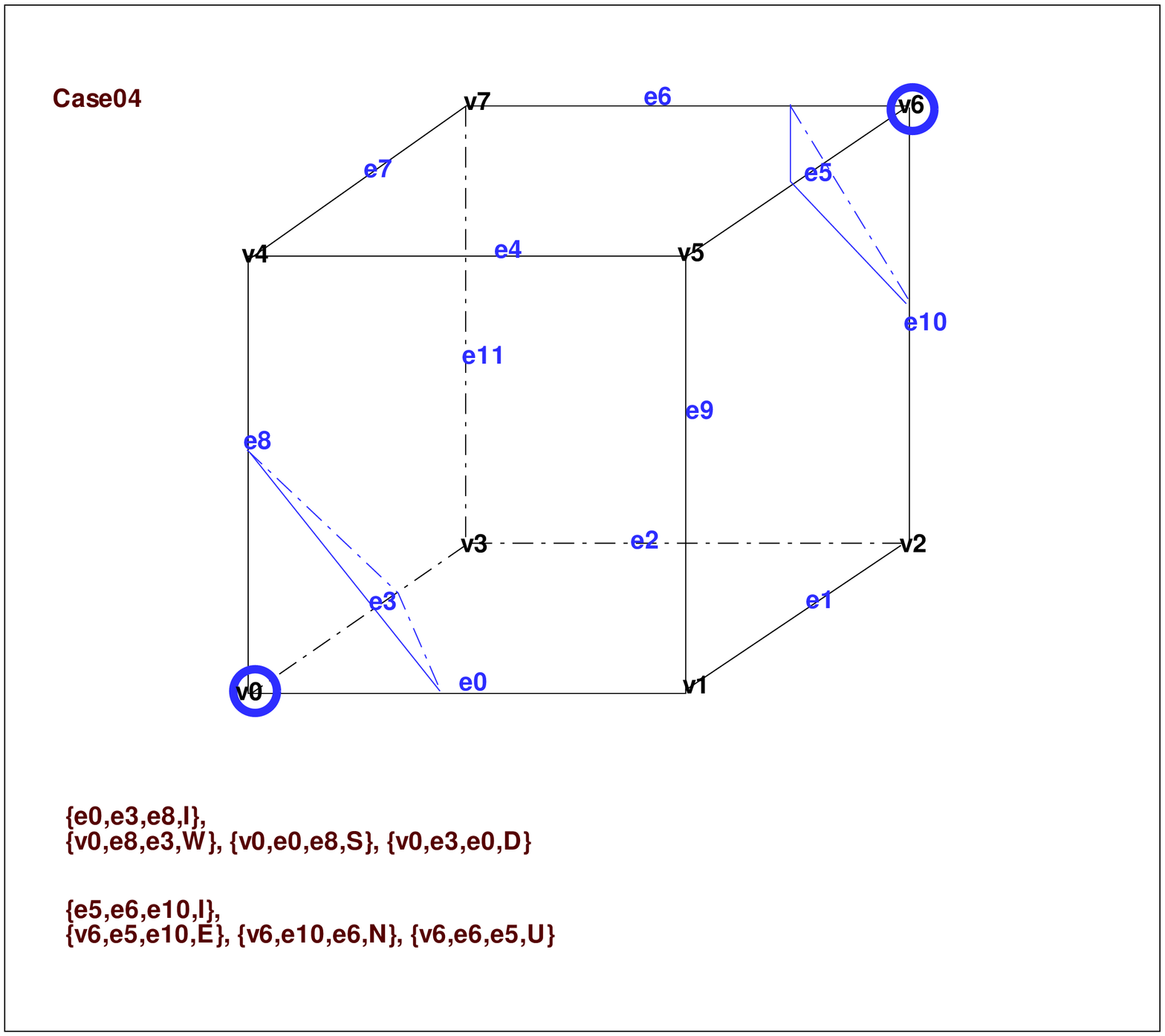} 
\end{figure}

\begin{figure}
\caption{\label{fig:Case05-Case06}Case05 and Case06}

\includegraphics[scale=0.3]{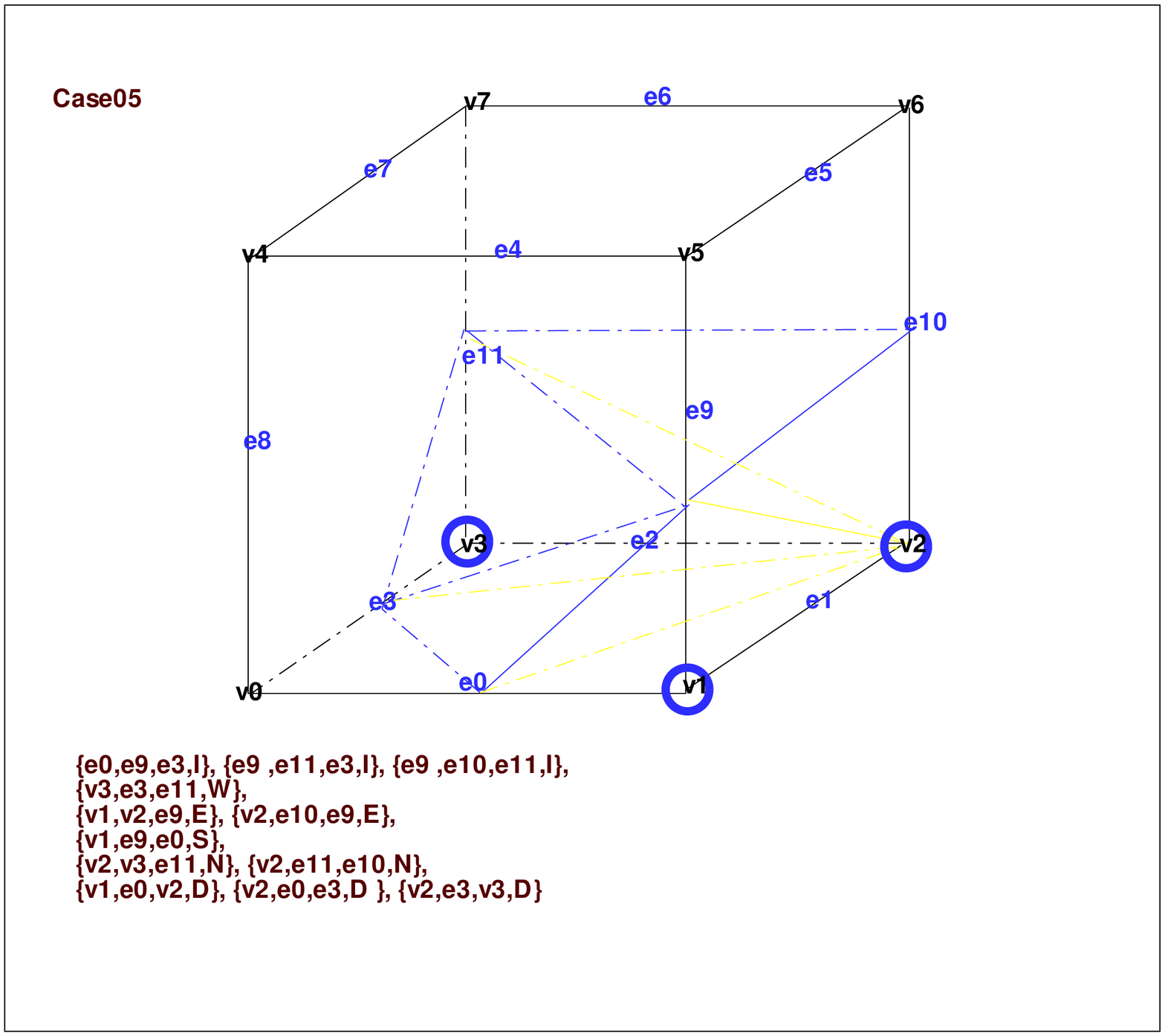} \includegraphics[scale=0.3]{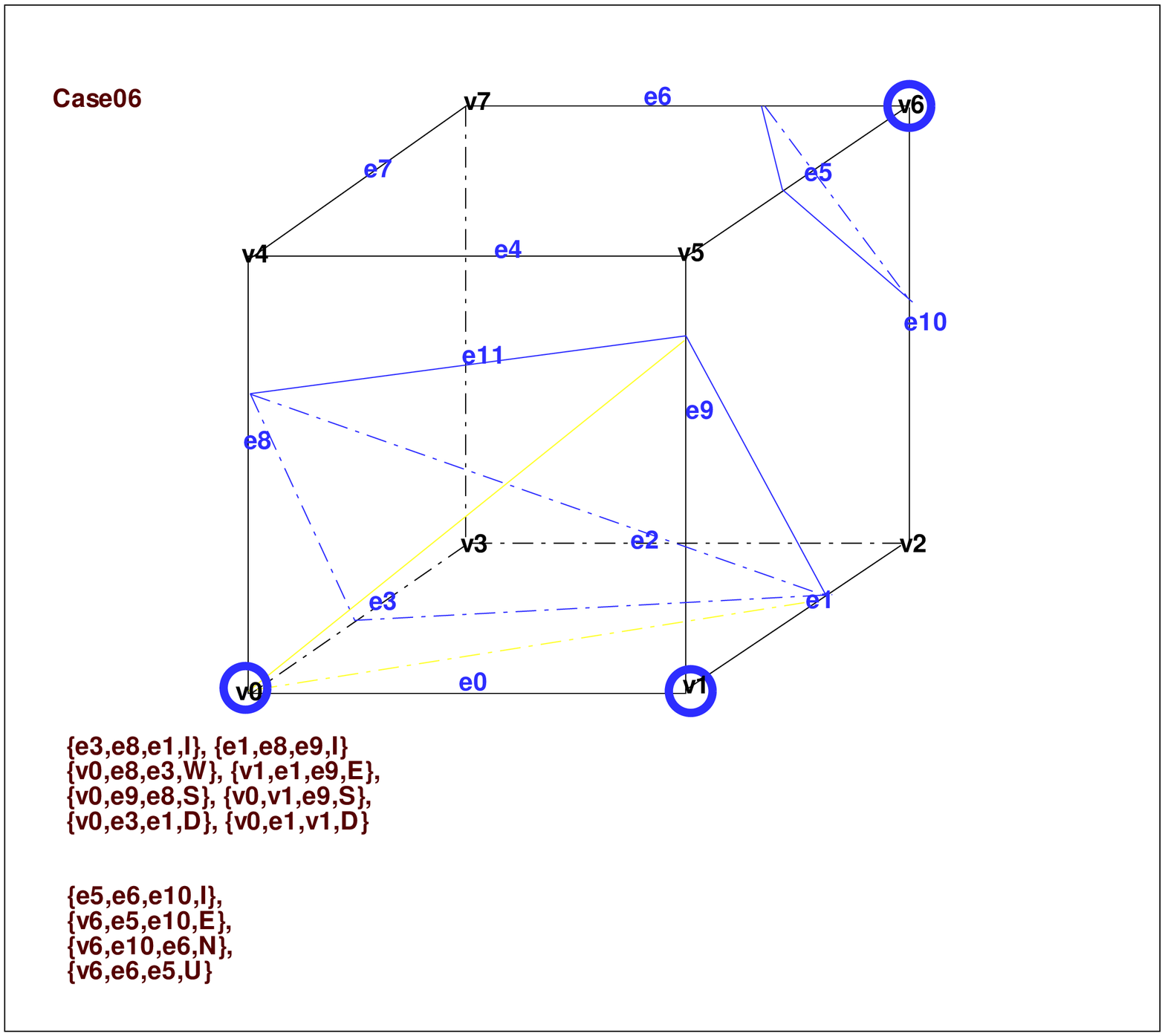} 
\end{figure}

\begin{figure}
\caption{\label{fig:Case07-Case08}Case07 and Case08}

\includegraphics[scale=0.3]{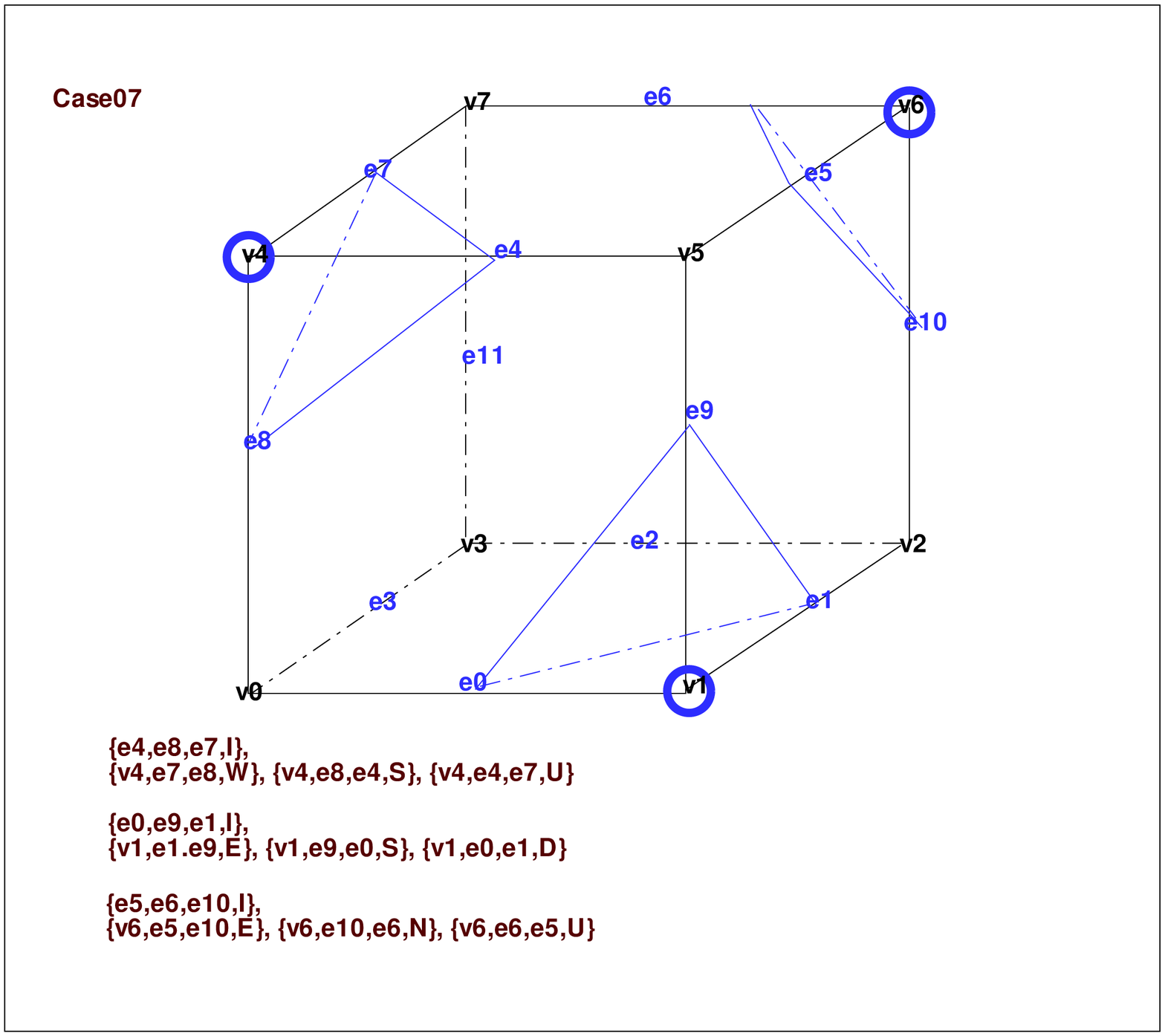} \includegraphics[scale=0.3]{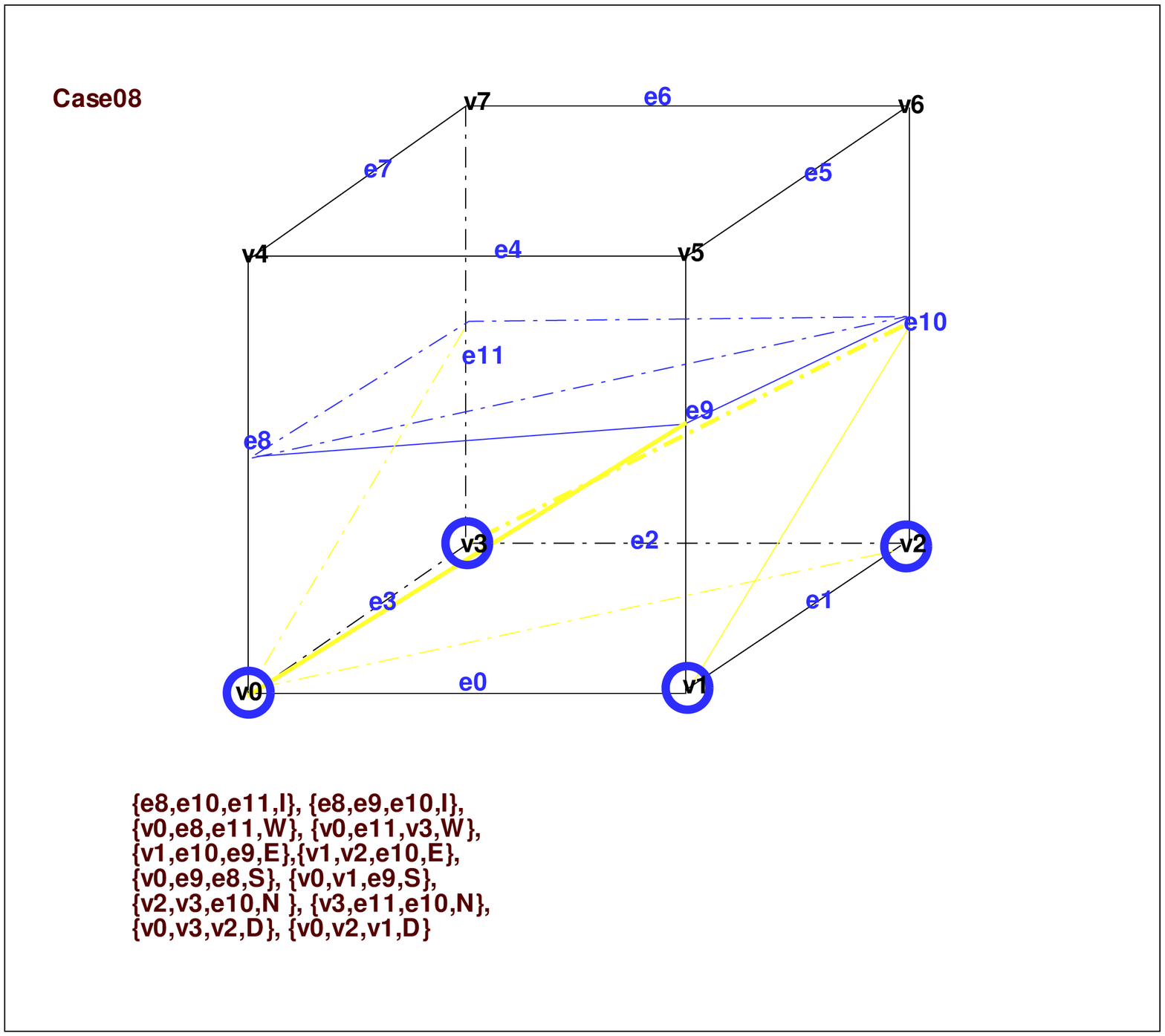} 
\end{figure}

\begin{figure}
\caption{\label{fig:Case09-Case10}Case09 and Case10}

\includegraphics[scale=0.3]{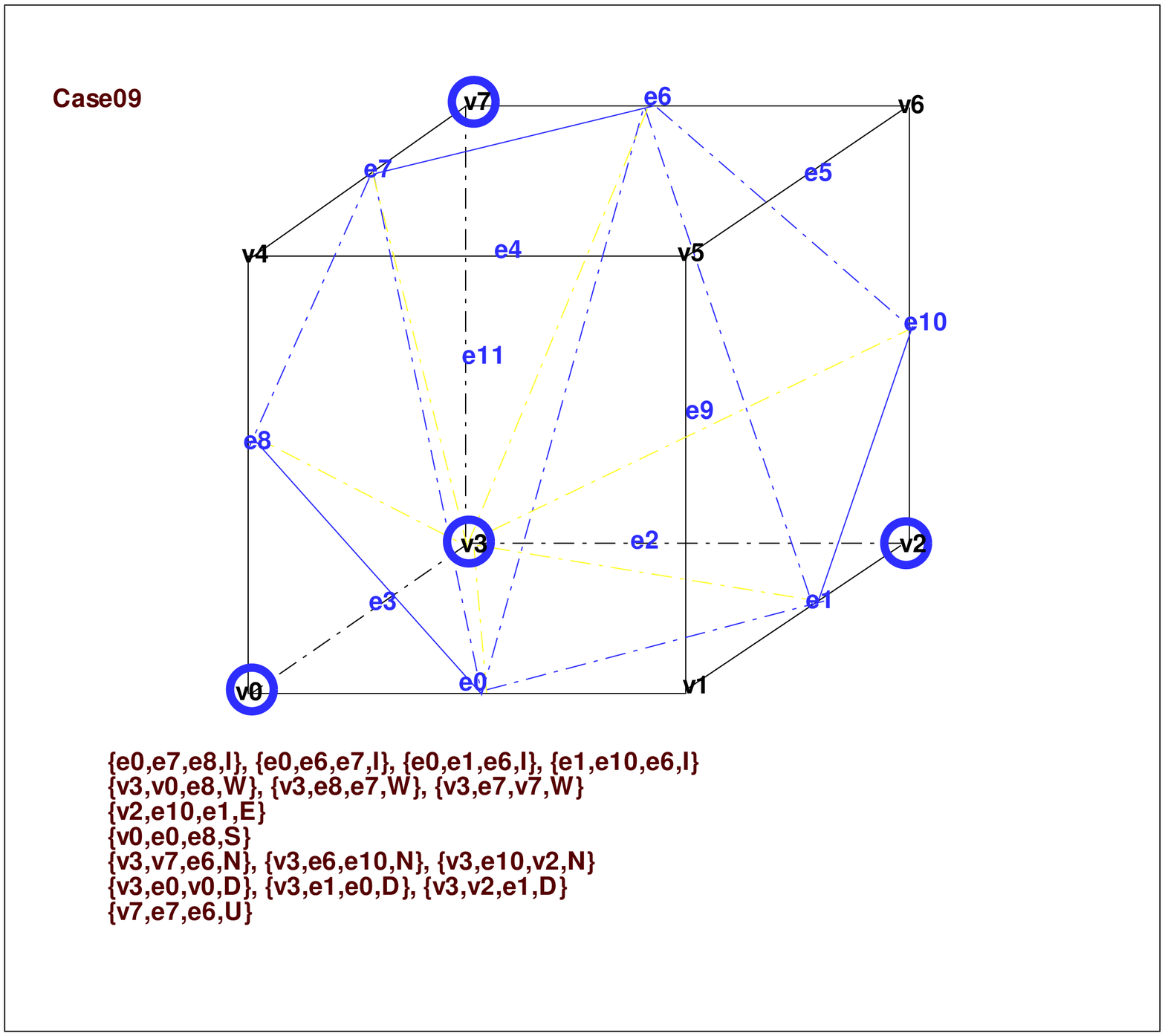} \includegraphics[scale=0.3]{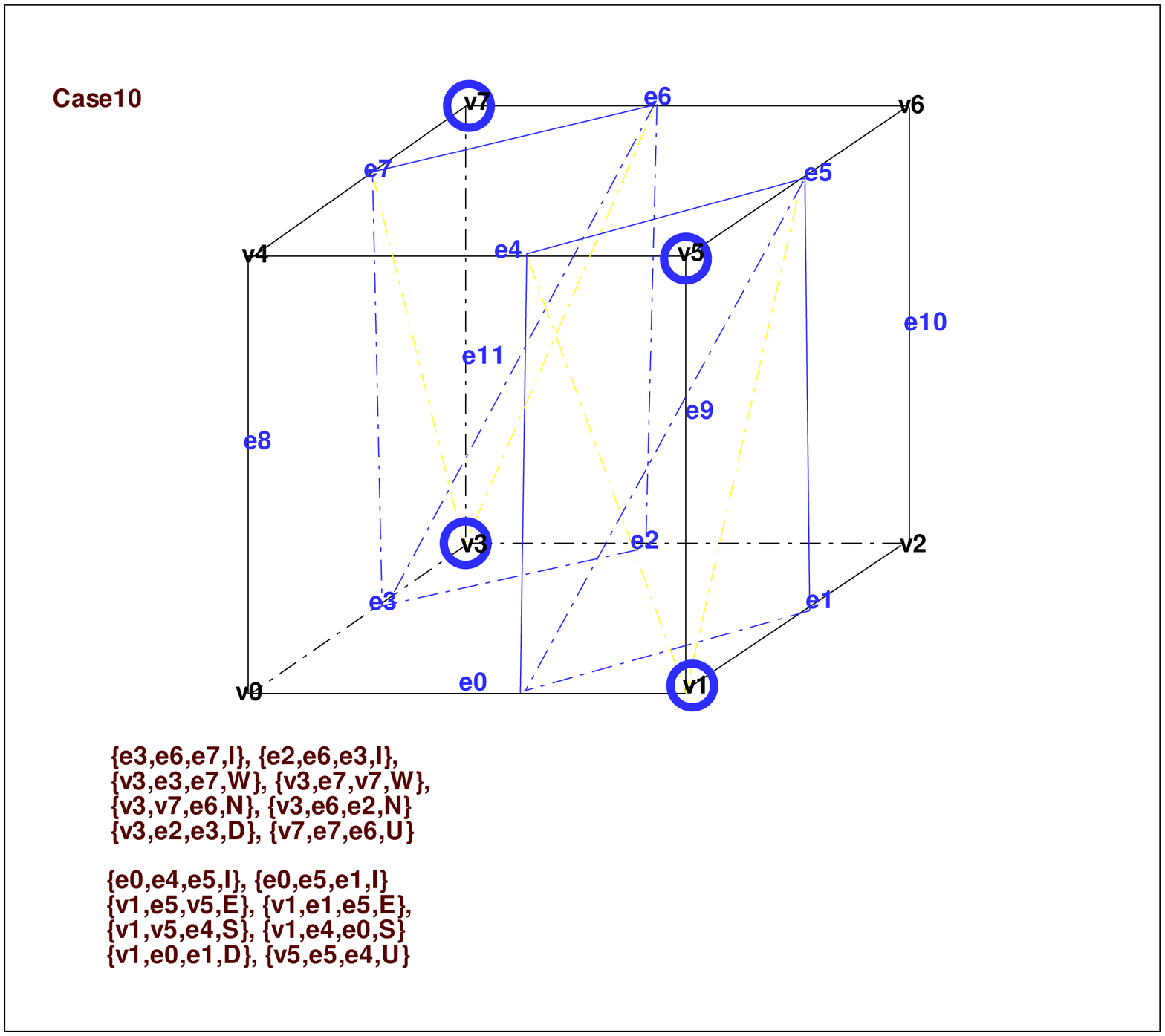} 
\end{figure}

\begin{figure}
\caption{\label{fig:Case11-Case12}Case11 and Case12}

\includegraphics[scale=0.3]{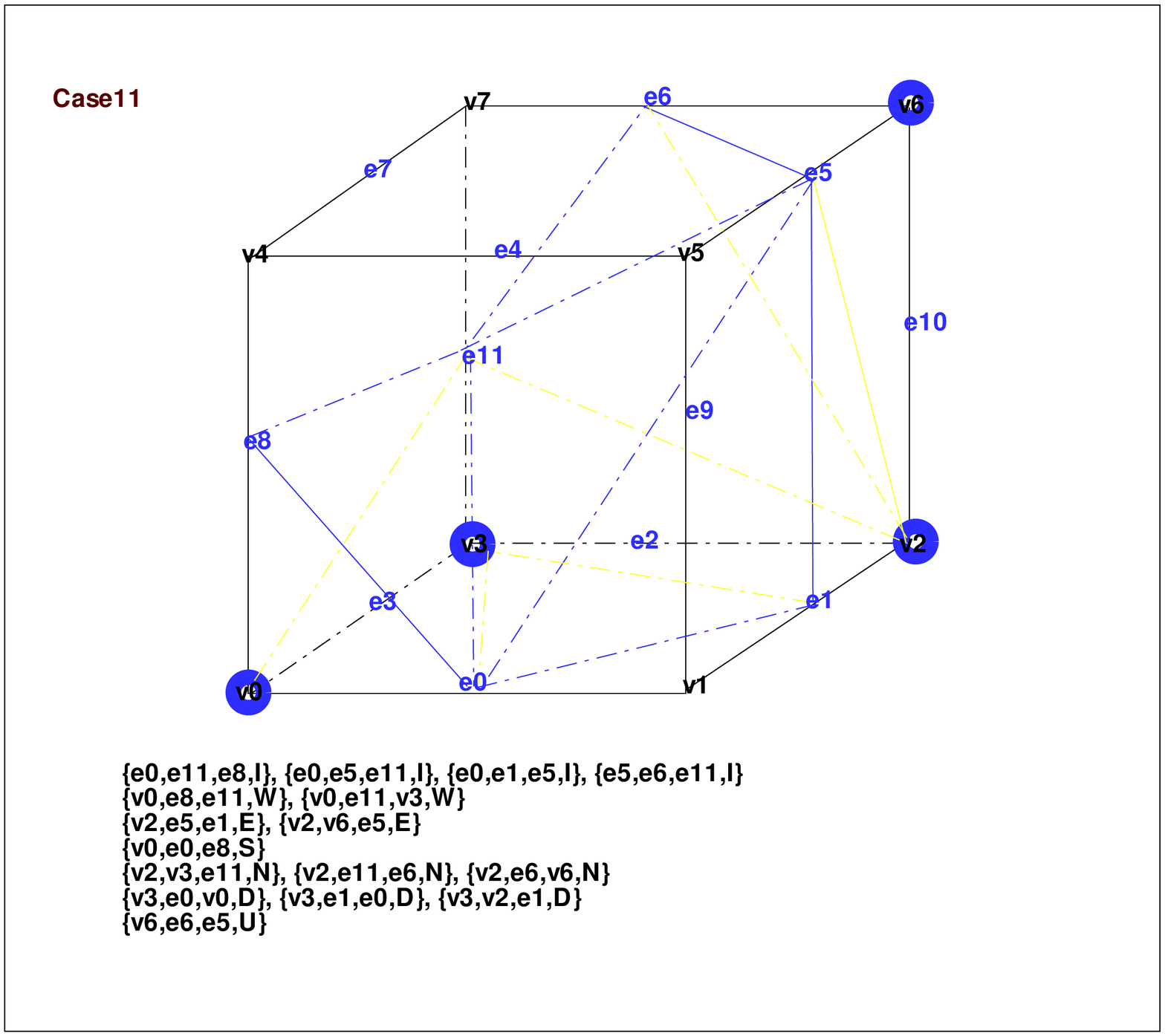} \includegraphics[scale=0.3]{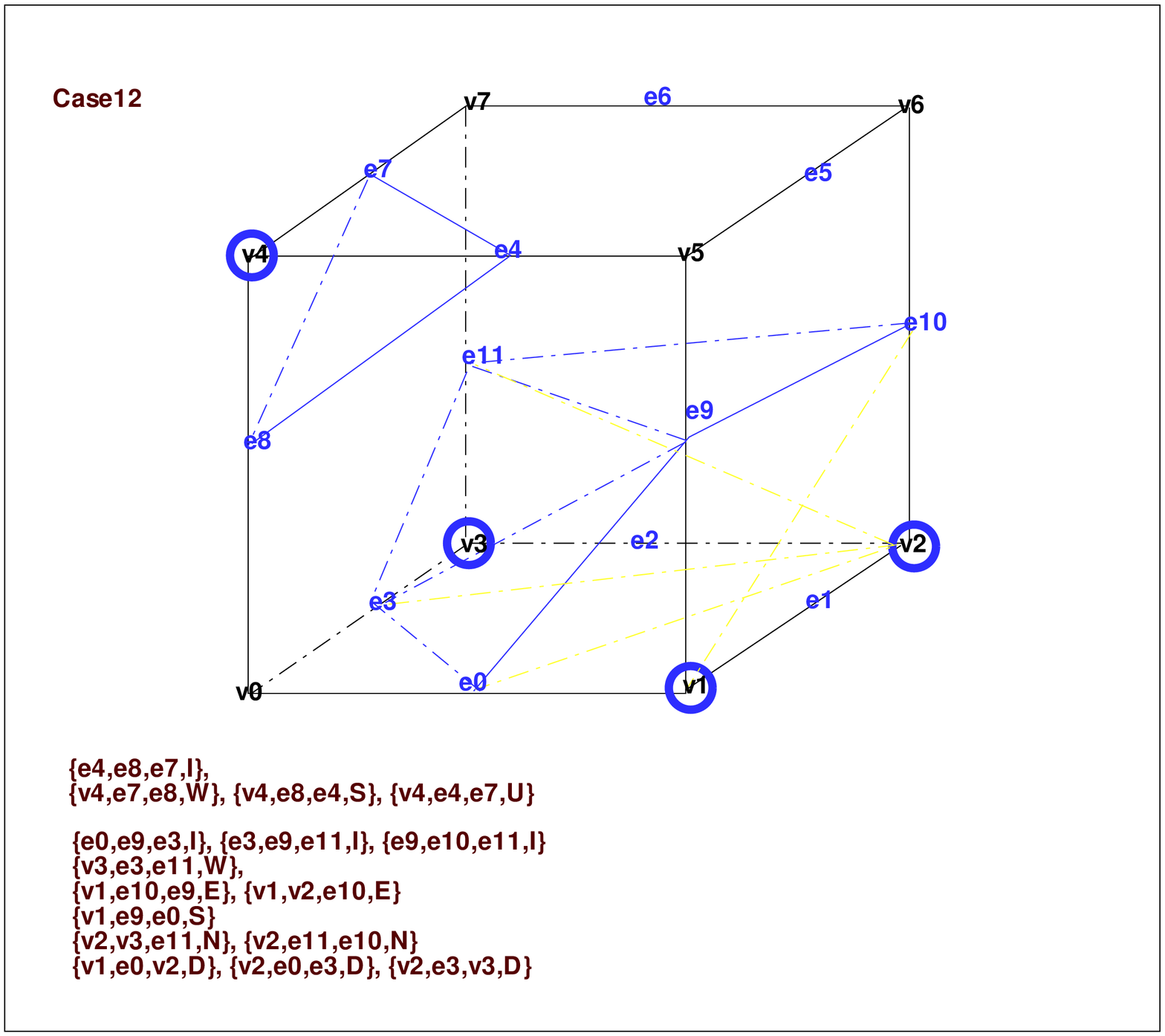} 
\end{figure}

\begin{figure}
\caption{\label{fig:Case13-Case14}Case13 and Case14}

\includegraphics[scale=0.3]{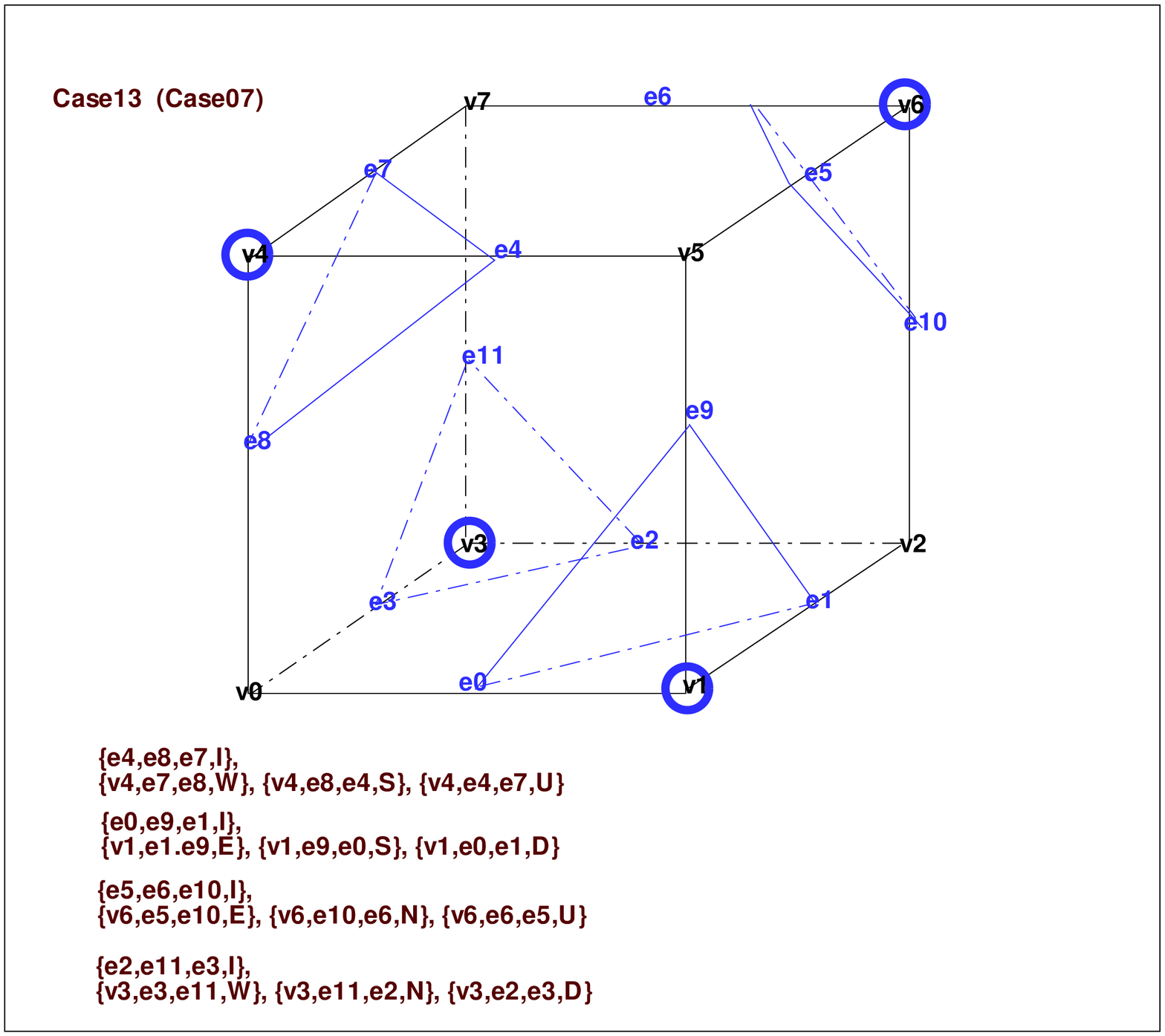} \includegraphics[scale=0.3]{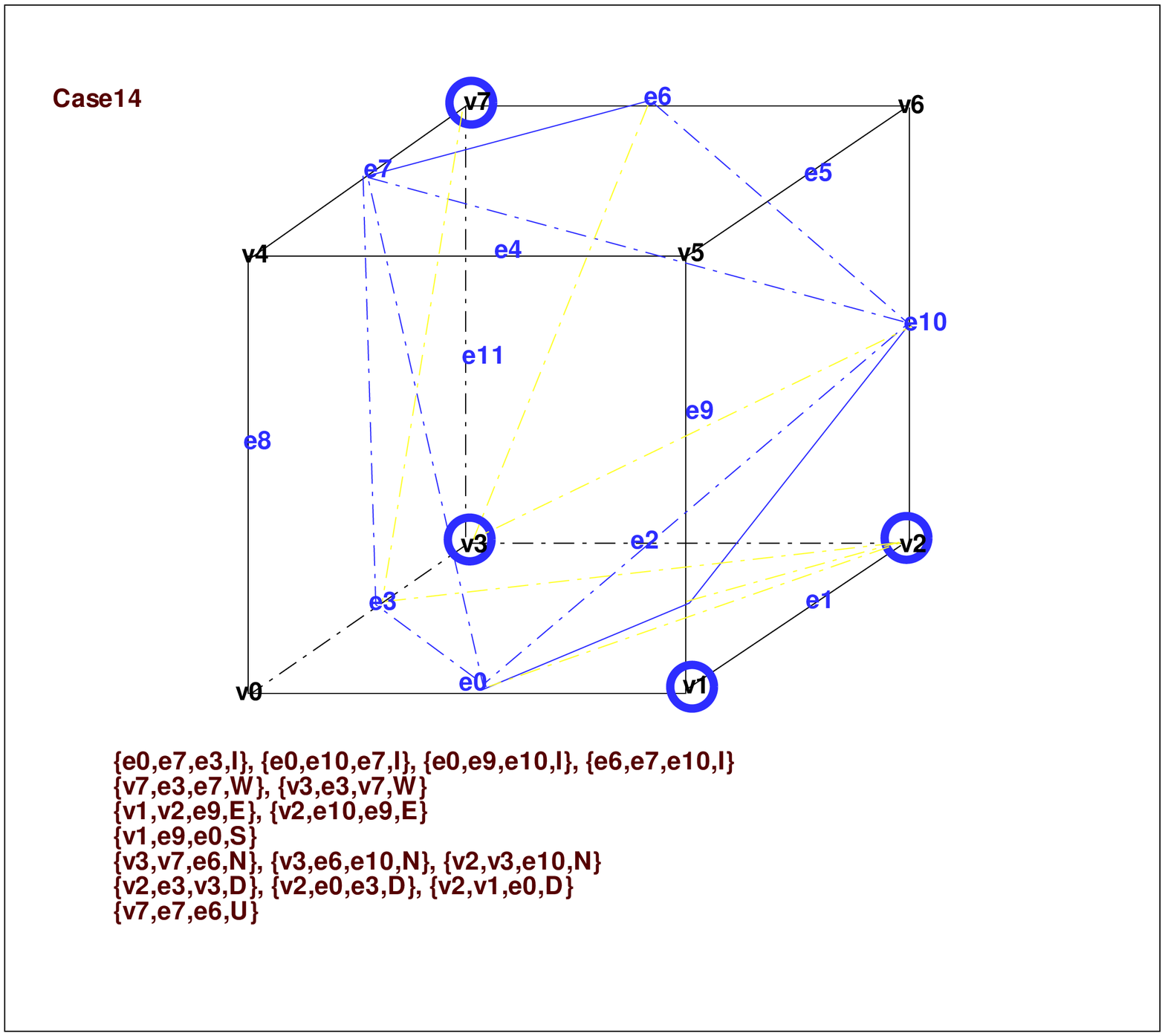} 
\end{figure}

\begin{figure}
\caption{\label{fig:Case15-Case16}Case15 and Case16}

\includegraphics[scale=0.3]{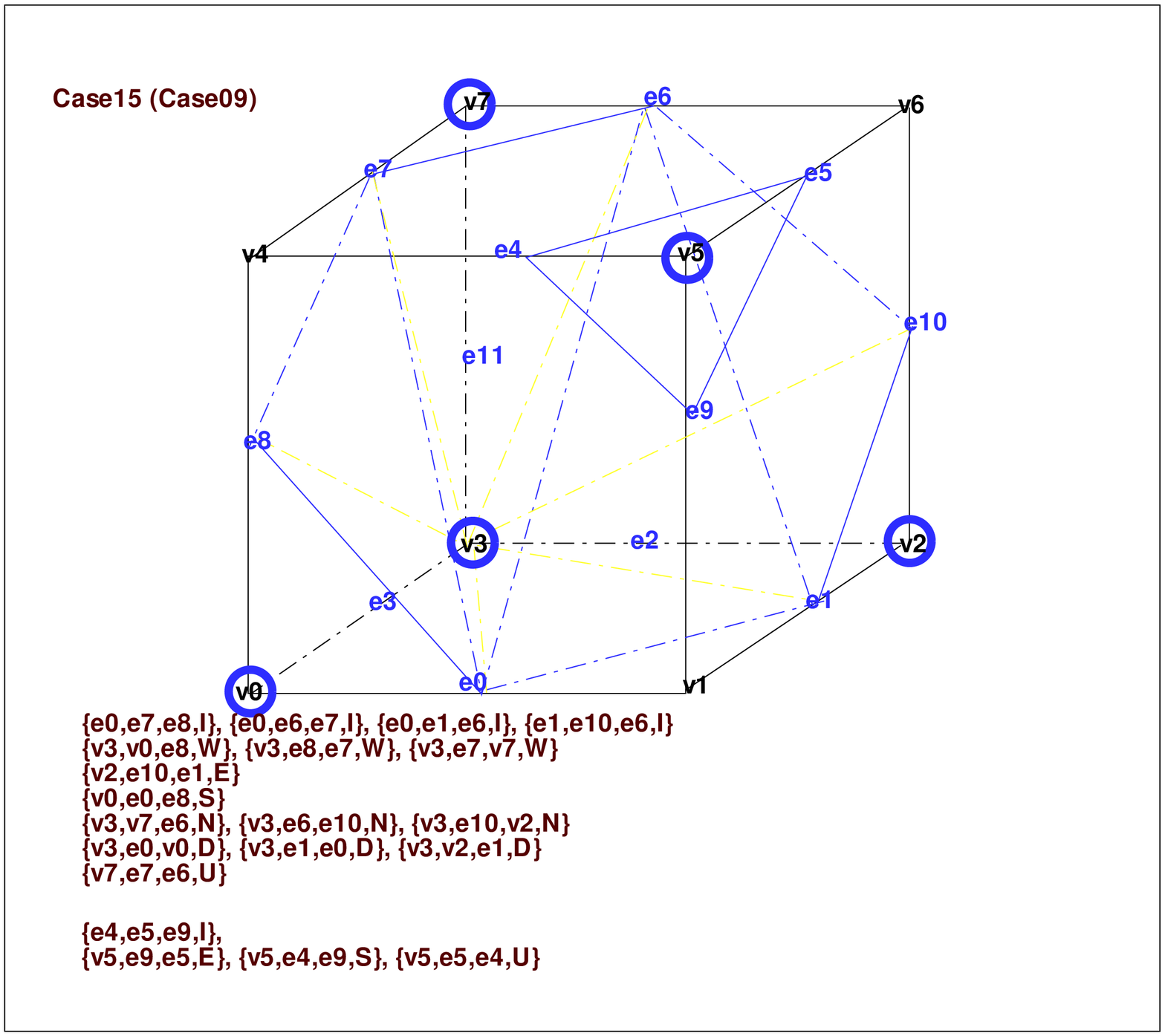} \includegraphics[scale=0.3]{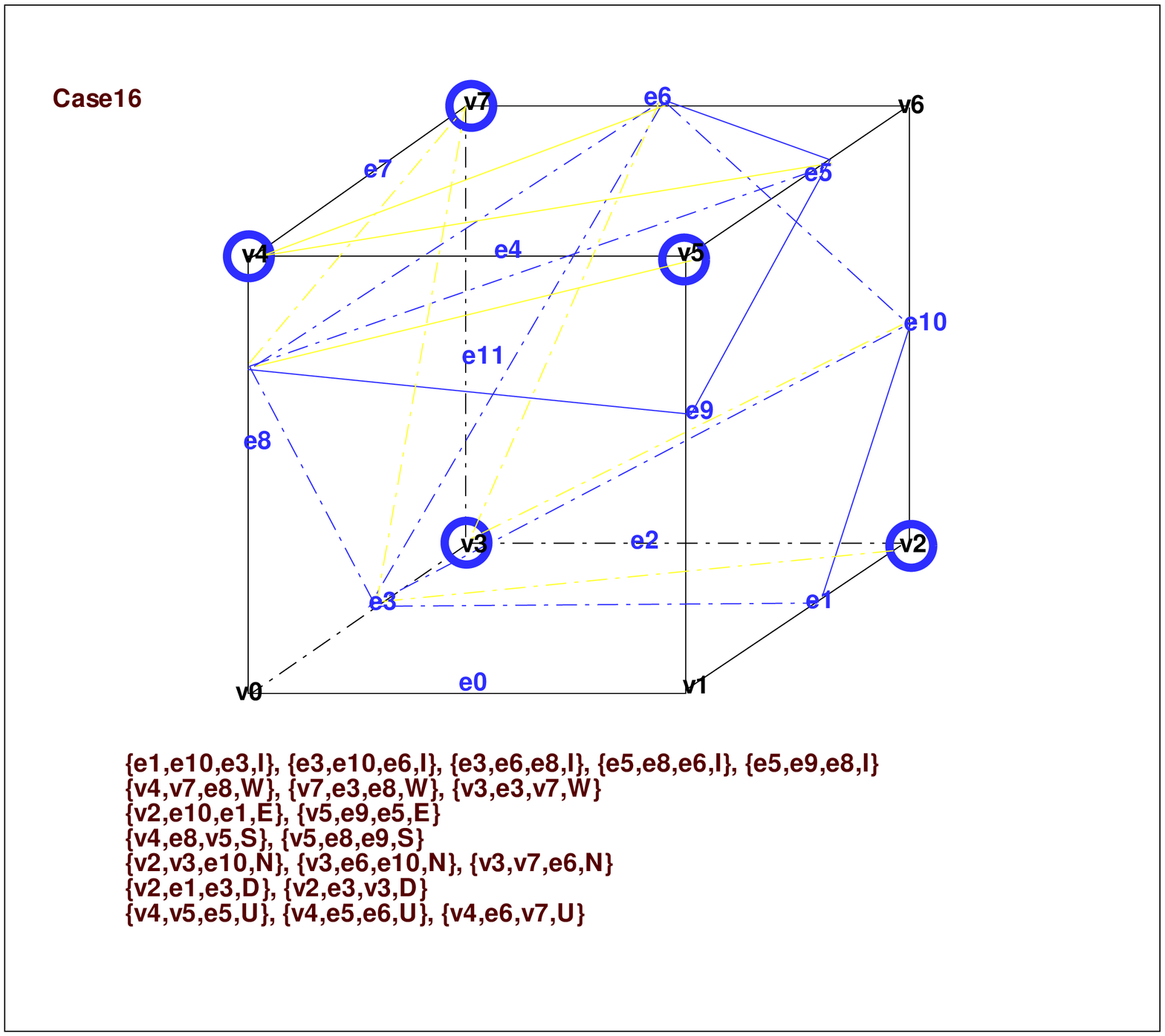} 
\end{figure}

\begin{figure}
\caption{\label{fig:Case17-Case18}Case17 and Case18}

\includegraphics[scale=0.3]{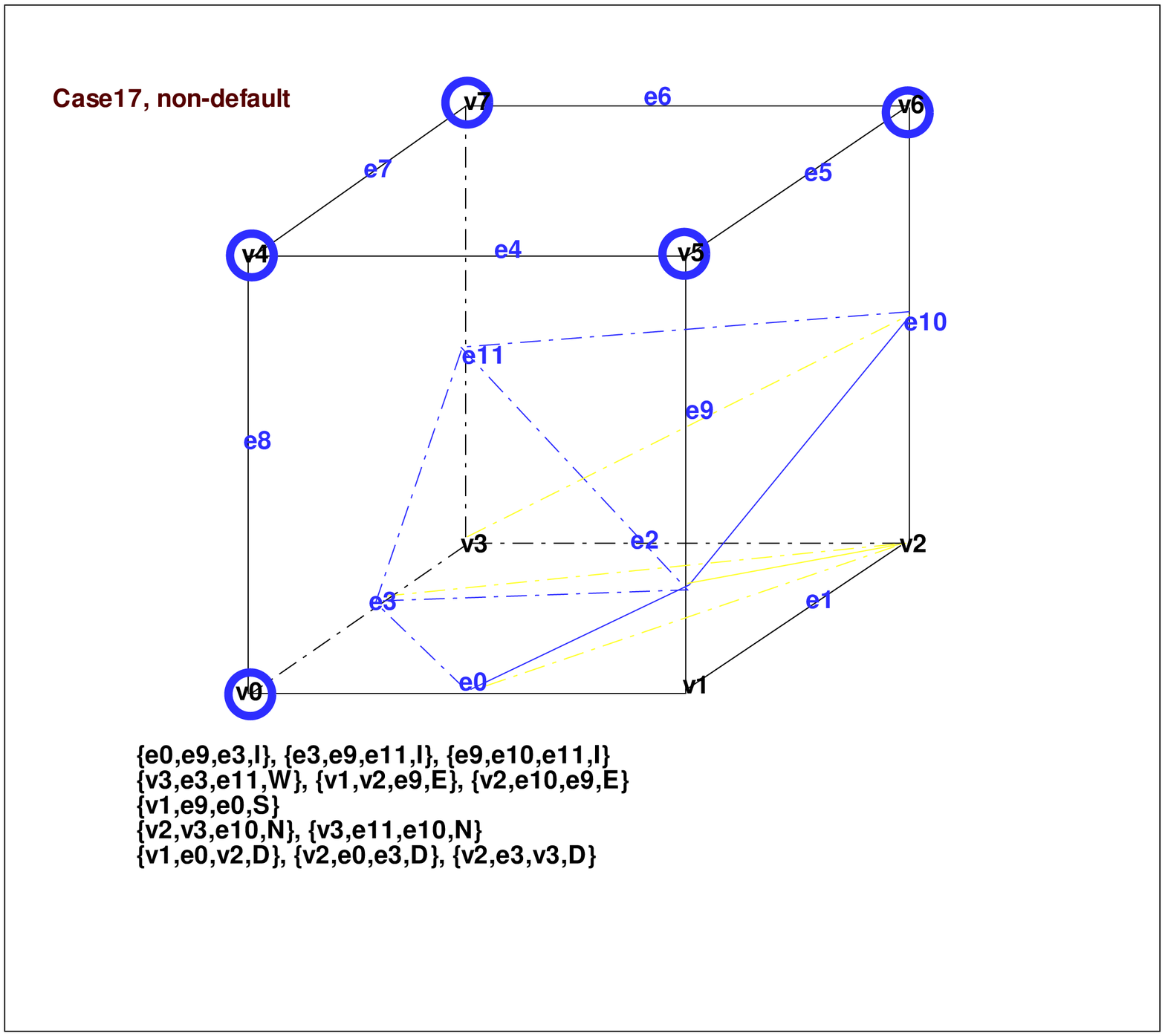} \includegraphics[scale=0.3]{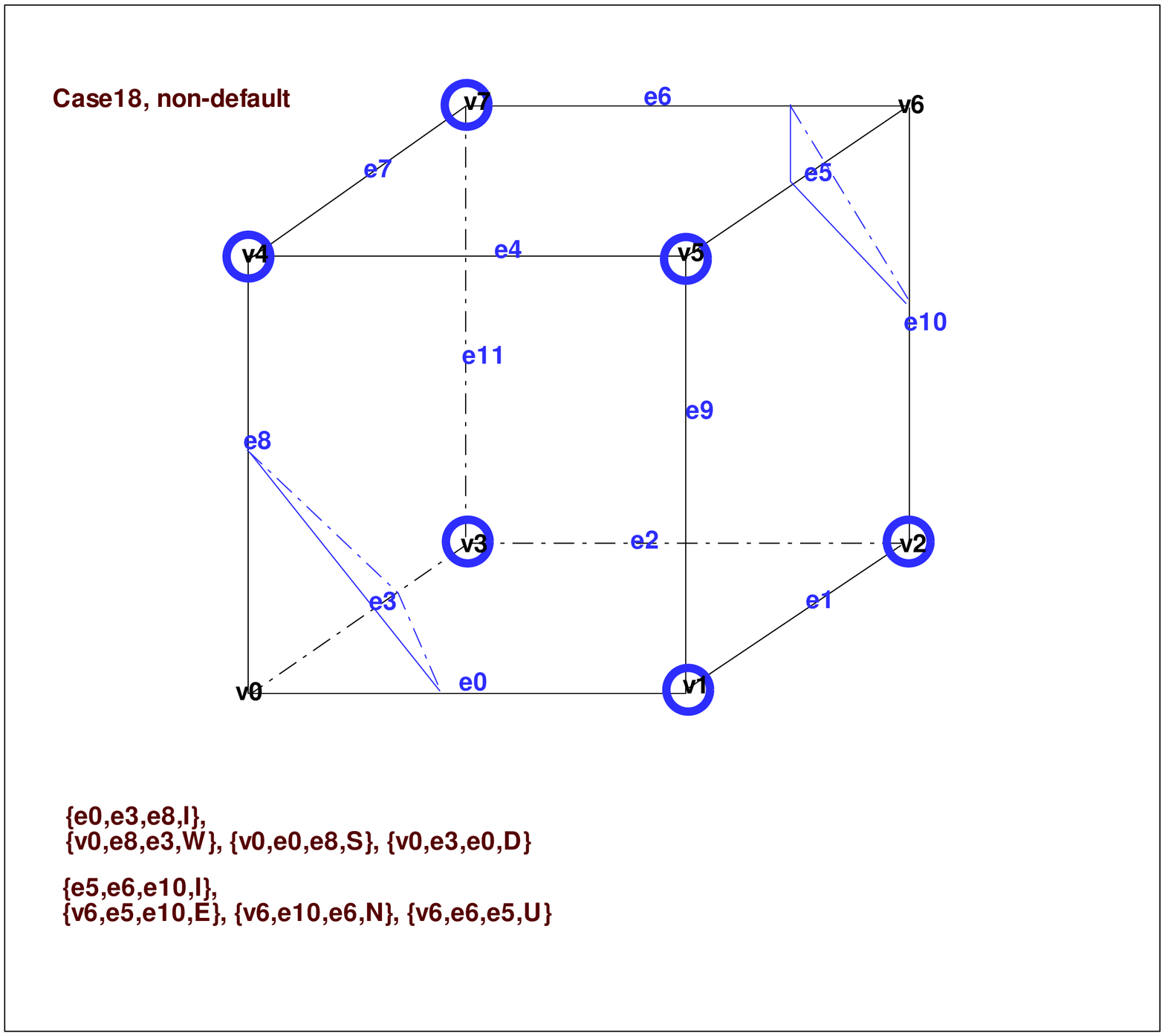} 
\end{figure}

\begin{figure}
\caption{\label{fig:Case19-Case20}Case19 and Case20}

\includegraphics[scale=0.3]{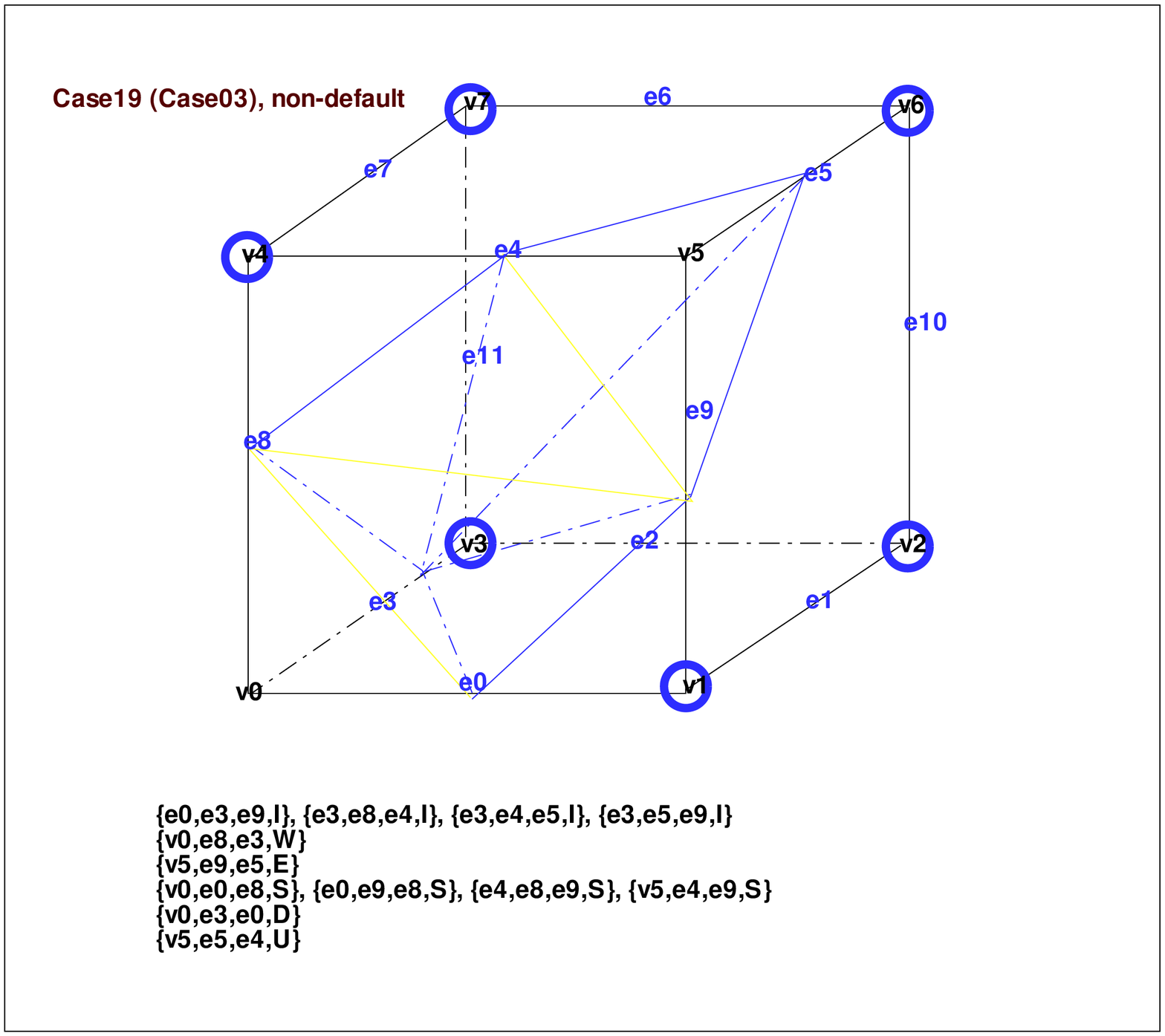} \includegraphics[scale=0.3]{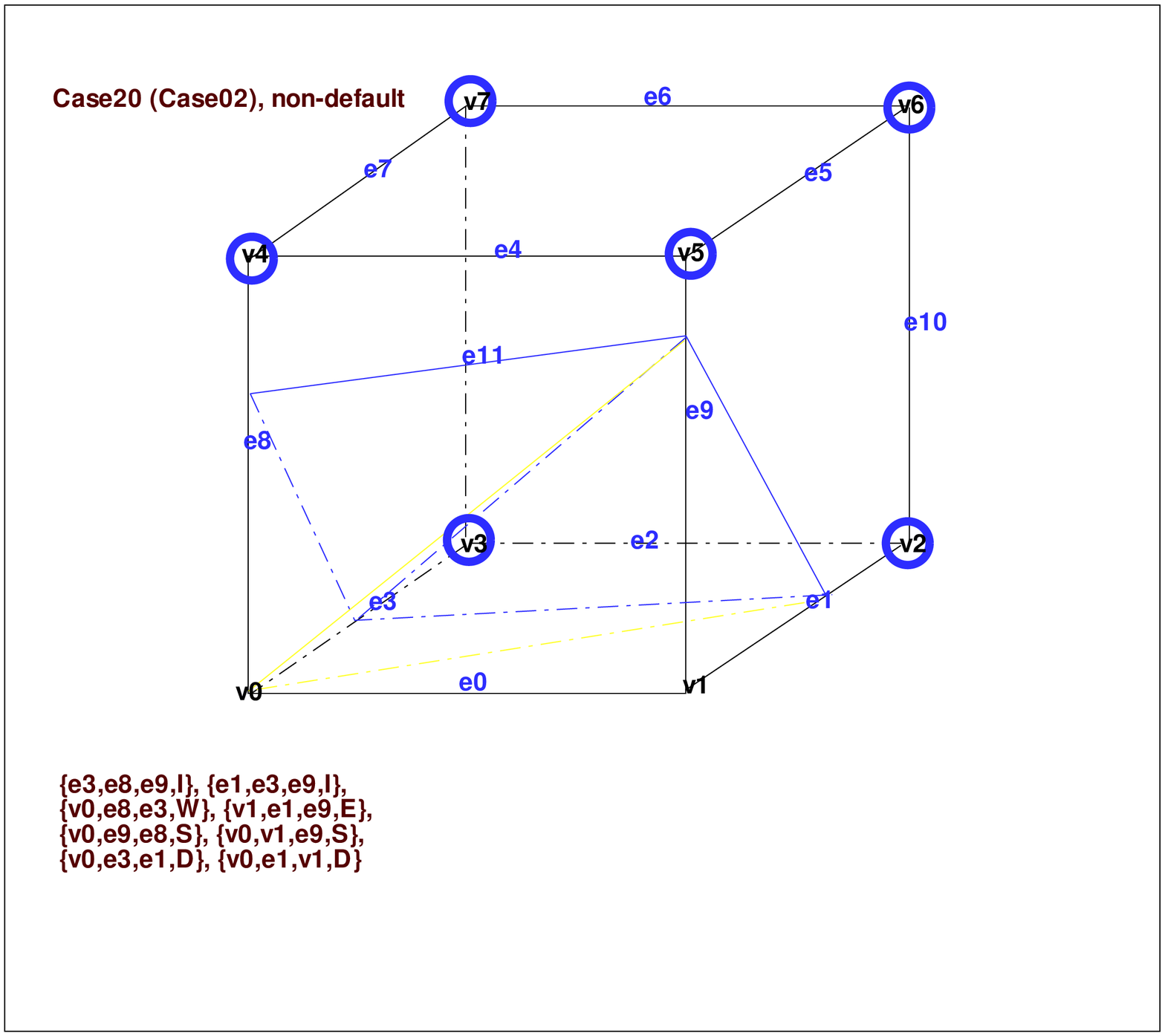} 
\end{figure}

\begin{figure}
\caption{\label{fig:Case21-Case22}Case21 and Case22}

\includegraphics[scale=0.3]{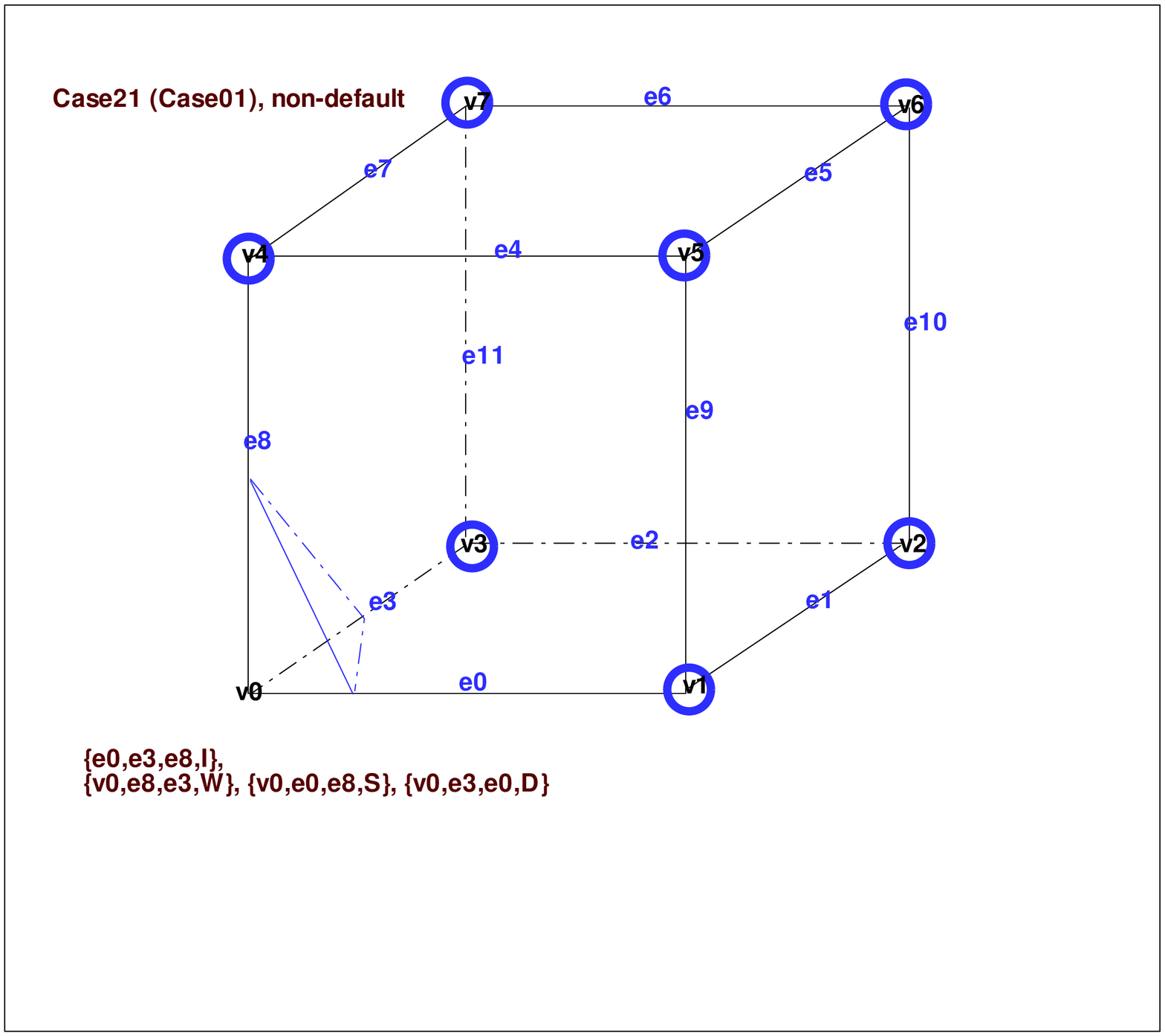} \includegraphics[scale=0.3]{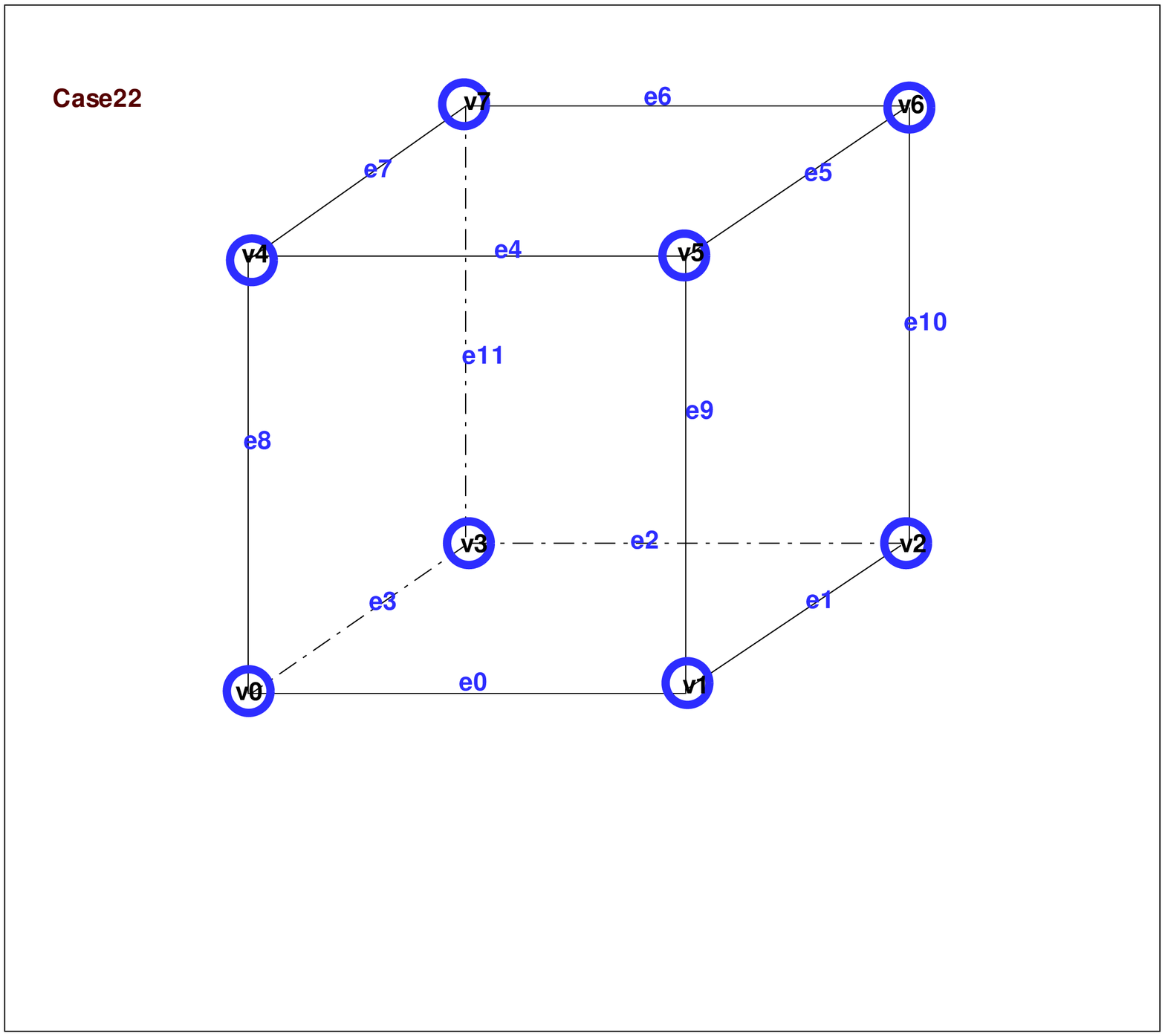} 
\end{figure}

For cubes with only two components, there are $2^{8}=256$ different
unique configurations. Using rotation symmetry only, these $256$
cases can be rotated into $23$ unique cases \cite{NewmanYi2006}
as shown in Figure \ref{fig:Case00}-\ref{fig:Case21-Case22}.

For each case, we have shown in the corresponding figures the triangulation
for the interface following \cite{NewmanYi2006} and the triangulations
for the six faces. Each triangle of the triangulations consist of
four letters: the first three letters are either the vertex index
or the edge index, and the fourth letter denotes the position of the
triangle: $I$ for triangle on the constructed interface, $W$ for
the triangle on the West face, and similarly for the other letters.

For example, for case $01$ in Figure \ref{fig:Case01-Case02}, the
triangulations list consists of

\begin{minipage}[t]{1\columnwidth}%
\{e0,e3,e8,I\}, \{v0,e8,e3,W\}, \{v0,e0,e8,S\}, \{v0,e3,e0,D\}%
\end{minipage}

For this case, there is only one triangle $\left\{ e0,e3,e8,I\right\} $
on the interface, one triangle $\left\{ v0,e8,e3,W\right\} $ on the
West face, one triangle $\left\{ v0,e0,e8,S\right\} $ on the South
face, and one triangle $\left\{ v0,e3,e0,D\right\} $ on the Down
face. The four triangles make up the closed surface enclosing the
domain corresponding to component $1$. Thus, to calculate the volume
for component $1$ for case $01$, we only need to use the surface
integration in equation \ref{eq:volume-calculation-for-Cartesian-Coordinate}.
Note that, the surface consists of four triangles, thus the surface
integration is a sum of the integrations on the four triangles, which
can be calculated easily.

For another example, say case $03$ in Figure \ref{fig:Case03-Case04},
the domain for component $1$ consists of two disconnected parts.
In the triangulation list, we first list the triangulation for the
first part, then the triangulations for the second part. In this way,
we could easily calculate other values besides the total volume for
each components, such as the number of the connected interfaces, the
averaged interface normal, interface center, and interface area for
each connected interface. If necessary, we can calculate separately
the partial volumes for the disconnected part.

For case $01-16$, we list the triangulations for component $1$.
For case $17-21$, we list the triangulations for component $0$ instead
of component $1$ since there are less triangles due to cube faces
triangulations. When we list the triangles for one connected domain,
we first list the triangles on the interface, then triangles on the
West, East, South, North, Down, and Upper faces.

After we have the triangulation list of the enclosing surface, the
volume in equation \ref{eq:volume-calculation-for-Cartesian-Coordinate}
can be calculated using \begin{eqnarray*}
Volume & = & \frac{\int_{\partial\Omega}\overrightarrow{n}\cdot\overrightarrow{f}ds}{3}\\
 & = & \frac{1}{3}\sum_{T_{i}}\int_{T_{i}}\overrightarrow{n}\cdot\overrightarrow{f}ds\end{eqnarray*}
 where $T_{i}$ is one of the triangle in the triangulation list.
Thus, the surface integration has been divided into a sum of the integration
on the triangles, which could be solved exactly using standard quadrature
rule.

It is also easy to calculate the interface center, interface normal,
interface area using the following formula:\begin{eqnarray*}
\vec{X}_{intfcCenter} & = & \int_{intfc}xds\\
\vec{N}_{intfcNormal} & = & \int_{intfc}\vec{n}ds\\
A_{intfcArea} & = & \int_{intfc}ds\end{eqnarray*}
 We simply use the triangulation of the interface only and then calculate
these surface integrations by the sum of the triangle integration.

\section{Examples}

For plane interfaces, numerical simulations show that our algorithm
give exact results for the total volumes and surface areas. In the
following, we use a sphere interface to do a mesh convergence study.

The computational domain is $[0,3]\times[0,3]\times[0,3]$. The interface
position is a sphere, given as\[
\sqrt{(x-1.5)^{2}+(y-1.5)^{2}+(z-1.5)^{2}}=1.\]
Table \ref{tab:3D-Mesh-Convergence-Study} shows the mesh convergence
study for the sphere volume and the surface area. From the table,
we can see that the method is second order accurate. 

\begin{table}
\caption{3D \label{tab:3D-Mesh-Convergence-Study}Mesh Convergence Study with
exact sphere volume $V=4.18879020479$ and sphere surface area $S=12.5663706144$}

\centering{}\begin{tabular}{|c|c|c|c|c|}
\hline 
Mesh Size & Volume & Error & Area & Error\tabularnewline
\hline
\hline 
10x10x10 & 4.00416024 & 0.18462996479 & 12.27248842 & 0.2938821944\tabularnewline
\hline 
20x20x20 & 4.13823907 & 0.05055113479 & 12.48614607 & 0.0802245444\tabularnewline
\hline 
40x40x40 & 4.17534124 & 0.01344896479 & 12.54510172 & 0.0212688944\tabularnewline
\hline 
80x80x80 & 4.18536131 & 0.00342889479 & 12.56094478 & 0.0054258344\tabularnewline
\hline
\end{tabular}
\end{table}

\section{Conclusion}

In this paper, we used the divergence theorem to calculate the partial
cell volumes. Compared with the tetrahedralization approach, it is
robust, efficient, and much easier for user to write computer code.
It is also straight forward to use this method to deal with more complicated
configurations \cite{Nielson2003}. 

The current method can be extended to Cylindrical coordinate, Spherical
coordinate. These will be documented in detail in future papers \cite{WangPolarCylindrical2013,WangSpherical2013}
since it is not easy to make them right.

\section*{Appendix Surface Triangles List}

We list the triangulation lists for all $23$ cases which can be easily
copied into C code. Note that the triangulation lists of Case $01-16$
are for component $1$ while the triangulation lists of Case $17-21$
are for component $0$.
\begin{itemize}
\item Case 00: N/A;
\item Case 01: \{ \{e0,e3,e8,I\}, \{v0,e8,e3,W\}, \{v0,e0,e8,S\}, \{v0,e3,e0,D\}
\};
\item Case 02: \{ \{e3,e8,e9,I\}, \{e1,e3,e9,I\}, \{v0,e8,e3,W\}, \{v1,e1,e9,E\},
\{v0,e9,e8,S\}, \{v0,v1,e9,S\}, \{v0,e3,e1,D\}, \{v0,e1,v1,D\} \};
\item Case 03: \{ \{e0,e3,e8,I\}, \{v0,e8,e3,W\}, \{v0,e0,e8,S\}, \{v0,e3,e0,D\},
\{e4,e5,e9,I\}, \{v5,e9,e5,E\}, \{v5,e4,e9,S\}, \{v5,e5,e4,U\} \};
\item Case 04: \{ \{e0,e3,e8,I\}, \{v0,e8,e3,W\}, \{v0,e0,e8,S\}, \{v0,e3,e0,D\},
\{e5,e6,e10,I\}, \{v6,e5,e10,E\}, \{v6,e10,e6,N\}, \{v6,e6,e5,U\}
\};
\item Case 05: \{ \{e0,e9,e3,I\}, \{e9,e11,e3,I\}, \{e9,e10,e11,I\}, \{v3,e3,e11,W\},
\{v1,v2,e9,E\}, \{v2,e10,e9,E\}, \{v1,e9,e0,S\}, \{v2,v3,e11,N\},
\{v2,e11,e10,N\}, \{v1,e0,v2,D\}, \{v2,e0,e3,D\}, \{v2,e3,v3,D\} \};
\item Case 06: \{ \{e3,e8,e1,I\}, \{e1,e8,e9,I\}, \{v0,e8,e3,W\}, \{v1,e1,e9,E\},
\{v0,e9,e8,S\}, \{v0,v1,e9,S\}, \{v0,e3,e1,D\}, \{v0,e1,v1,D\}, \{e5,e6,e10,I\},
\{v6,e5,e10,E\}, \{v6,e10,e6,N\}, \{v6,e6,e5,U\} \};
\item Case 07: \{ \{e4,e8,e7,I\}, \{v4,e7,e8,W\}, \{v4,e8,e4,S\}, \{v4,e4,e7,U\},
\{e0,e9,e1,I\}, \{v1,e1,e9,E\}, \{v1,e9,e0,S\}, \{v1,e0,e1,D\}, \{e5,e6,e10,I\},
\{v6,e5,e10,E\}, \{v6,e10,e6,N\}, \{v6,e6,e5,U\} \};
\item Case 08: \{ \{e8,e10,e11,I\}, \{e8,e9,e10,I\}, \{v0,e8,e11,W\}, \{v0,e11,v3,W\},
\{v1,e10,e9,E\}, \{v1,v2,e10,E\}, \{v0,e9,e8,S\}, \{v0,v1,e9,S\},
\{v2,v3,e10,N\}, \{v3,e11,e10,N\}, \{v0,v3,v2,D\}, \{v0,v2,v1,D\}
\};
\item Case 09: \{ \{e0,e7,e8,I\}, \{e0,e6,e7,I\}, \{e0,e1,e6,I\}, \{e1,e10,e6,I\},
\{v3,v0,e8,W\}, \{v3,e8,e7,W\}, \{v3,e7,v7,W\}, \{v2,e10,e1,E\}, \{v0,e0,e8,S\},
\{v3,v7,e6,N\}, \{v3,e6,e10,N\}, \{v3,e10,v2,N\}, \{v3,e0,v0,D\},
\{v3,e1,e0,D\}, \{v3,v2,e1,D\}, \{v7,e7,e6,U\} \};
\item Case 10: \{ \{e3,e6,e7,I\}, \{e2,e6,e3,I\}, \{v3,e3,e7,W\}, \{v3,e7,v7,W\},
\{v3,v7,e6,N\}, \{v3,e6,e2,N\}, \{v3,e2,e3,D\}, \{v7,e7,e6,U\}, \{e0,e4,e5,I\},
\{e0,e5,e1,I\}, \{v1,e5,v5,E\}, \{v1,e1,e5,E\}, \{v1,v5,e4,S\}, \{v1,e4,e0,S\},
\{v1,e0,e1,D\}, \{v5,e5,e4,U\} \};
\item Case 11: \{ \{e0,e11,e8,I\}, \{e0,e5,e11,I\}, \{e0,e1,e5,I\}, \{e5,e6,e11,I\},
\{v0,e8,e11,W\}, \{v0,e11,v3,W\}, \{v2,e5,e1,E\}, \{v2,v6,e5,E\},
\{v0,e0,e8,S\}, \{v2,v3,e11,N\}, \{v2,e11,e6,N\}, \{v2,e6,v6,N\},
\{v3,e0,v0,D\}, \{v3,e1,e0,D\}, \{v3,v2,e1,D\}, \{v6,e6,e5,U\} \};
\item Case 12: \{ \{e4,e8,e7,I\}, \{v4,e7,e8,W\}, \{v4,e8,e4,S\}, \{v4,e4,e7,U\},
\{e0,e9,e3,I\}, \{e3,e9,e11,I\}, \{e9,e10,e11,I\}, \{v3,e3,e11,W\},
\{v1,e10,e9,E\}, \{v1,v2,e10,E\}, \{v1,e9,e0,S\}, \{v2,v3,e11,N\},
\{v2,e11,e10,N\}, \{v1,e0,v2,D\}, \{v2,e0,e3,D\}, \{v2,e3,v3,D\} \};
\item Case 13: \{ \{e4,e8,e7,I\}, \{v4,e7,e8,W\}, \{v4,e8,e4,S\}, \{v4,e4,e7,U\},
\{e0,e9,e1,I\}, \{v1,e1,e9,E\}, \{v1,e9,e0,S\}, \{v1,e0,e1,D\}, \{e5,e6,e10,I\},
\{v6,e5,e10,E\}, \{v6,e10,e6,N\}, \{v6,e6,e5,U\}, \{e2,e11,e3,I\},
\{v3,e3,e11,W\}, \{v3,e11,e2,N\}, \{v3,e2,e3,D\} \};
\item Case 14: \{ \{e0,e7,e3,I\}, \{e0,e10,e7,I\}, \{e0,e9,e10,I\}, \{e6,e7,e10,I\},
\{v7,e3,e7,W\}, \{v3,e3,e7,W\}, \{v1,v2,e9,E\}, \{v2,e10,e9,E\}, \{v1,e9,e0,S\},
\{v3,v7,e6,N\}, \{v3,e6,e10,N\}, \{v2,v3,e10,N\}, \{v2,e3,v3,D\},
\{v2,e0,e3,D\}, \{v2,v1,e0,D\}, \{v7,e7,e6,U\} \};
\item Case 15: \{ \{e0,e7,e8,I\}, \{e0,e6,e7,I\}, \{e0,e1,e6,I\}, \{e1,e10,e6,I\},
\{v3,v0,e8,W\}, \{v3,e8,e7,W\}, \{v3,e7,v7,W\}, \{v2,e10,e1,E\}, \{v0,e0,e8,S\},
\{v3,v7,e6,N\}, \{v3,e6,e10,N\}, \{v3,e10,v2,N\}, \{v3,e0,v0,D\},
\{v3,e1,e0,D\}, \{v3,v2,e1,D\}, \{v7,e7,e6,U\}, \{e4,e5,e9,I\}, \{v5,e9,e5,E\},
\{v5,e4,e9,S\}, \{v5,e5,e4,U\} \};
\item Case 16: \{ \{e1,e10,e3,I\}, \{e3,e10,e6,I\}, \{e3,e6,e8,I\}, \{e5,e8,e6,I\},
\{e5,e9,e8,I\}, \{v4,v7,e8,W\}, \{v7,e3,e8,W\}, \{v3,e3,v7,W\}, \{v2,e10,e1,E\},
\{v5,e9,e5,E\}, \{v4,e8,v5,S\}, \{v5,e8,e9,S\}, \{v2,v3,e10,N\}, \{v3,e6,e10,N\},
\{v3,v7,e6,N\}, \{v2,e1,e3,D\}, \{v2,e3,v3,D\}, \{v4,v5,e5,U\}, \{v4,e5,e6,U\},
\{v4,e6,v7,U\} \};
\item Case 17: \{ \{e0,e9,e3,I\}, \{e3,e9,e11,I\}, \{e9,e10,e11,I\}, \{v3,e3,e11,W\},
\{v1,v2,e9,E\}, \{v2,e10,e9,E\}, \{v1,e9,e0,S\}, \{v2,v3,e10,N\},
\{v3,e11,e10,N\}, \{v1,e0,v2,D\}, \{v2,e0,e3,D\}, \{v2,e3,v3,D\} \};
\item Case 18: \{ \{e0,e3,e8,I\}, \{v0,e8,e3,W\}, \{v0,e0,e8,S\}, \{v0,e3,e0,D\},
\{e5,e6,e10,I\}, \{v6,e5,e10,E\}, \{v6,e10,e6,N\}, \{v6,e6,e5,U\}
\};
\item Case 19: \{ \{e0,e3,e9,I\}, \{e3,e8,e4,I\}, \{e3,e4,e5,I\}, \{e3,e5,e9,I\},
\{v0,e8,e3,W\}, \{v5,e9,e5,E\}, \{v0,e0,e8,S\}, \{e0,e9,e8,S\}, \{e4,e8,e9,S\},
\{v5,e4,e9,S\}, \{v0,e3,e0,D\}, \{v5,e5,e4,U\} \};
\item Case 20: \{ \{e3,e8,e9,I\}, \{e1,e3,e9,I\}, \{v0,e8,e3,W\}, \{v1,e1,e9,E\},
\{v0,e9,e8,S\}, \{v0,v1,e9,S\}, \{v0,e3,e1,D\}, \{v0,e1,v1,D\} \};
\item Case 21: \{ \{e0,e3,e8,I\}, \{v0,e8,e3,W\}, \{v0,e0,e8,S\}, \{v0,e3,e0,D\}
\};
\item Case 22: N/A.
\end{itemize}

\section*{Appendix Rotation List }

For completeness, we list the 24 unique rotation of a cube \cite{MathsCubeRotation}
in the following:
\begin{itemize}
\item \{0,1,2,3,4,5,6,7\}, // self
\item \{4,5,1,0,7,6,2,3\}, // opposite face: x 90
\item \{7,6,5,4,3,2,1,0\}, // opposite face: x 180
\item \{3,2,6,7,0,1,5,4\}, // opposite face: x 270
\item \{4,0,3,7,5,1,2,6\}, // opposite face: y 90
\item \{5,4,7,6,1,0,3,2\}, // opposite face: y 180
\item \{1,5,6,2,0,4,7,3\}, // opposite face: y 270
\item \{3,0,1,2,7,4,5,6\}, // opposite face: z 90
\item \{2,3,0,1,6,7,4,5\}, // opposite face: z 180
\item \{1,2,3,0,5,6,7,4\}, // opposite face: z 270
\item \{0,4,5,1,3,7,6,2\}, // opposite vertices: v0-v6
\item \{0,3,7,4,1,2,6,5\}, // opposite vertices: v0-v6
\item \{2,1,5,6,3,0,4,7\}, // opposite vertices: v1-v7
\item \{5,1,0,4,6,2,3,7\}, // opposite vertices: v1-v7
\item \{5,6,2,1,4,7,3,0\}, // opposite vertices: v2-v4
\item \{7,3,2,6,4,0,1,5\}, // opposite vertices: v2-v4
\item \{2,6,7,3,1,5,4,0\}, // opposite vertices: v3-v5
\item \{7,4,0,3,6,5,1,2\}, // opposite vertices: v3-v5
\item \{1,0,4,5,2,3,7,6\}, // opposite lines: e0-e6
\item \{3,7,4,0,2,6,5,1\}, // opposite lines: e3-e5
\item \{6,7,3,2,5,4,0,1\}, // opposite lines: e2-e4
\item \{6,2,1,5,7,3,0,4\}, // opposite lines: e1-e7
\item \{4,7,6,5,0,3,2,1\}, // opposite lines: e8-e10
\item \{6,5,4,7,2,1,0,3\} // opposite lines: e9-e11
\end{itemize}
\bibliographystyle{plain} \bibliographystyle{plain}
\bibliography{../refs}

\end{document}